\documentclass[11pt]{article}
\usepackage{amssymb, amsmath}
\usepackage{color}
\usepackage{graphics}
\usepackage{graphicx}
\usepackage{amssymb}
\usepackage{epsf}
\usepackage{verbatim}
\usepackage{color}
\usepackage{hyperref}

\usepackage[normalem]{ulem}

\newtheorem{theorem}{Theorem}
\newtheorem{lemma}{Lemma}

\newtheorem{definition}[lemma]{Definition}
\newtheorem{proposition}[lemma]{Proposition}

\newtheorem{example}[lemma]{Example}
\newtheorem{remark}[lemma]{Remark}

\newcommand{\vone}{\vec{1}}
\newcommand{\vs}{\vec{s}}
\newcommand{\vsst}{\vec{s}^{\, *}}

\newcommand{\vp}{\vec{p}}
\newcommand{\vth}{\vec{th}}

\newcommand{\cC}{{\mathcal C}}

\newcommand{\cP}{{\mathcal P}}

\newcommand{\cU}{{\mathcal U}}
\newcommand{\eps}{\varepsilon}

\definecolor{OliveGreen}{cmyk}{0.64,0,0.95,0.40}
\definecolor{gold}{rgb}{0.85,.66,0}

\begin{document}

\title{Lengths of attractors and transients in neuronal networks with random connectivities}
\author{ Winfried Just\footnote{Department of Mathematics, Ohio University, Athens, OH 45701, USA,  mathjust@gmail.com}, Sungwoo Ahn\footnote{Department of Mathematical Sciences, Indiana University-Purdue University Indianapolis, IN 46202, USA, ahnmath@gmail.com}}

\maketitle

\bigskip

\begin{abstract}
We study how the dynamics of a class of  discrete dynamical system models for neuronal networks depends on the connectivity of the network. Specifically, we assume that the network is an Erd\H{o}s-R\'{enyi} random graph and analytically derive scaling laws for the average lengths of the attractors and transients under certain restrictions on the intrinsic parameters of the neurons, that is, their refractory periods and firing thresholds.  In contrast to earlier results that were reported
in~\cite{TAWJ}, here we focus on the connection probabilities near the phase transition where the most complex dynamics is expected to occur.
\end{abstract}

\bigskip

\section{Introduction}\label{sec:intro}

We study generic behavior of a class of discrete dynamic system models of neuronal networks.  This class was designed to model the phenomenon of so-called \emph{dynamic clustering,} where time appears to be partitioned into consecutive episodes during which certain groups of neurons fire together, while other neurons are quiescent and membership in the respective groups may change from episode to episode~\cite{ahn_terman_2010, bazhenov_stopfer_2001}.  This phenomenon has been observed in a number of actual neuronal tissues~\cite{laurent1996, stopfer1997}.  In~\cite{TAWJ} it was proved that dynamic clustering occurs in a broad class of ODE models with a certain architecture, and that, moreover, discrete models as described here will reliably predict, for a large region of the state space of these ODE models, which neurons fire during a given episode.

In order to gain biological insights from this type of models, one needs to understand how the dynamics of the network depends on the network connectivity, which is modeled as a directed graph~$D$. Some provable restrictions on the possible network dynamics for certain classes of digraphs were derived in~\cite{thesis,AJ}
and Section~6.4 of~ \cite{bookchapter}. For most neuronal networks in actual organisms though, the precise connectivity is not known, and it is important to study the expected dynamics under the assumption that~$D$ is randomly drawn from some probability distribution.  In~\cite{JAT} the case of Erd\H{o}s-Renyi random digraphs was investigated, and two phase transitions were characterized.  If the other parameters of the model are set to $p^* = th^* = 1$ (see below for definitions), these occur when the mean degrees are $\sim \ln n$ and $\approx 3$ respectively. However, as simulation studies such as those in~\cite{ahn_terman_2010} suggest, the most interesting dynamics are expected when the mean degree is~$\approx 1$.  The purpose of the present paper is to derive characterizations of the phase transition in the expected dynamics that occur for the latter type of connectivities.

\section{Basic definitions}\label{sec:defs}

Our notation will be mostly standard.  The number of elements of a finite set~$A$ will be denoted by~$|A|$.  All digraphs considered here are implicitly assumed to be loop-free and have vertex set
$[n] = \{1, \ldots , n\}$; the arc set of a digraph~$D$ will be denoted by~$A_D$. A directed path from~$i_0$ to~$i_\ell$ of length~$\ell$ will be identified with a sequence $(i_0, i_1, \ldots , i_\ell)$ of pairwise distinct nodes such that $<i_j, i_{j+1}> \, \in A_D$ for all $j < \ell$.  Similarly, a directed cycle of length~$\ell$ will be identified with a sequence $(i_0, i_1, \ldots , i_\ell)$ of  nodes such that $<i_j, i_{j+1}> \, \in A_D$ for all $j < \ell$, the nodes $i_0, \ldots, i_{\ell-1}$ are pairwise distinct, and $i_0 = i_\ell$. We call a directed path~$Pt$  \emph{straight} if it has the property that for any two vertices $i,j$ in~$Pt$ there exists exactly one directed path in~$D$, either from~$i$ to~$j$ or \emph{vice versa.}  Note that this unique path must be a part of~$Pt$.

The set of all vertices~$j$ that can be reached from a given vertex~$i$ via a directed path in~$D$ will be denoted by~$DC_D(i)$, or $DC(i)$ if~$D$ is implied by the context, and called the \emph{downstream component of~$i$ (in~$D$);} the dual object  will be denoted by~$UC_D(i)$, or $UC(i)$  and called the \emph{upstream component of~$i$ (in~$D$).} \emph{The strongly connected component of~$i$ (in~$D$)} is the set
$SC_D(i) = SC(i)  = DC(i) \cap UC(i)$. The set of vertices in the  largest (giant) strongly connected component of~$D$ will be denoted by~$GC$ or~$GC_D$. Similarly, $UG, DG$ will denote the corresponding upstream and downstream components.
We will sometimes  use the abbreviations~$UC(W) = \bigcup_{i \in W} UC(i)$ \emph{etc.}  for a subset~$W \subseteq [n]$.
 Moreover, to reduce clutter in our notation, we will somewhat informally also use $UC(i), GC$ \emph{etc.} as symbols for the subdigraphs of~$D$ that are induced by these sets of vertices. We hope that this will not cause confusion as it should always be clear from the context which meaning is intended.

 The maximum length of a directed path in~$D$ will be denoted by~$L_{\max}$,  the maximum length of a straight directed path in~$D$ will be denoted by~$L_{\max}^s$, and the maximum length of a directed path in~$UC(i)$ will be denoted by~$L_{\max}(i)$.

\begin{definition}\label{def:network:discrete_model}
A \emph{neuronal network} is a triple $N ~= ~<D, \vec{p},
\vec{th}>$, where $D = \ <[n], A_D>$ is a loop-free
digraph and $\vp = (p_1, \ldots , p_n)$ and $\vth = (th_1, \ldots ,
th_n)$ are vectors of positive integers. The state $\vs$ of the
system at time $t$ is a vector $\vec{s}(t) = (s_1(t), \ldots ,
s_n(t))$, where $s_i(t) \in \{0, 1, \ldots , p_i\}$ for all $i \in
[n]$. The state space of $N$ will be denoted by $St_N$. The updating function
of $N$ is defined as follows:
\begin{itemize}
\item[(R)] If $s_i(t) < p_i$, then $s_i(t+1) = s_i(t) + 1$.
\item[(F)] If $s_i(t) = p_i$ and there exist at least $th_i$ different $j \in [n]$ with  $s_j(t) = 0$ and
$<j, i> \ \in A_D$, then $s_i(t+1) = 0$.
\item[(E)] If $s_i(t) = p_i$ and there are fewer than $th_i$ different $j \in [n]$ with  $s_j(t) = 0$ and
$<j, i> \ \in A_D$, then $s_i(t+1) = p_i$.
\end{itemize}
\end{definition}

In terms of neuroscience applications, $s_i(t)=0$ represents the event that neuron~$i$ \emph{fires} at time~$t$.  The number~$p_i$ represents the \emph{refractory period} of neuron number~$i$, that is, the minimum number of time steps for which a neuron needs to remain in a resting state after firing until it can fire again. This is formalized by condition~(R) above. Once neuron~$i$ has reached the end of its refractory period, it will fire if it receives \emph{firing input} from a number of other neurons that is at least equal to its \emph{firing threshold} $th_i$ (condition~(F)), and remain at the end of its refractory period otherwise (condition~(E)).

Since the state space $St_N$ is finite, for each \emph{trajectory} $(\vec{s}(0), \vec{s}(1), \vec{s}(2) , \dots )$  there must exist times $0 \leq t_0 < t_1$ with
$\vec{s}(t_0) = \vec{s}(t_1)$.   The smallest time~$t_0$ for which this occurs will be denoted by~$\tau$.   It represents the
length of the \emph{transient,} that is, of the initial segment $(\vec{s}(0), \vec{s}(1), \dots, \vec{s}(\tau-1))$  of a trajectory  that consists of states that will be visited only once.  The set
$\{\vec{s}(\tau), \vec{s}(\tau +1), \dots, \vec{s}(2\tau-1) \}$ consists of all states  that are visited infinitely often. It is called the \emph{attractor} of the trajectory.  Its size will be denoted by~$\alpha$. Sometime it is more intuitive to think of the attractor as a sequence, and we will often write about the ``length'' of an attractor instead of its ``size.''  Each of our networks has exactly one attractor of length~1; its only element is the unique \emph{steady state}~$\vp$. Attractors of length~$\alpha > 1$ will be called \emph{periodic attractors}. \emph{The basin of attraction} of a given attractor is the set of all initial conditions from which it will be reached.  Note that in general both~$\alpha$ and~$\tau$ depend both on the particular network and on the initial state.  If the network, the initial state, or both are sampled from some distribution, then~$\alpha$ and~$\tau$ become r.v.s.

 Example~\ref{ex:Landau} of Appendix~A illustrates these concepts and may help readers to gain some familiarity with the dynamics of our networks.

Throughout this note we will assume that the connectivity~$D$ of a network $N ~= ~ < D, \vec{p}, \vec{th}>$ is randomly drawn from the distribution of~Erd\H{o}s-R\'{e}nyi random digraphs, such that each potential arc~$<i,j>$ is included in~$A_D$ with a given probability~$\pi(n)$,  the vectors $\vp$ or $\vth$ are randomly drawn from the uniform distributions of all vectors $\vp$, $\vth$ such that $p_* \leq \min \vp \leq \max \vp \leq p^*$ and $th_* \leq \min \vth \leq \max \vth \leq th^*$, where
$p_*, p^*, th_*, th^*$ are parameters that do not depend on~$n$, and initial states $\vs(0)$ are randomly drawn from the uniform distribution of all possible initial states in the resulting network.

We are interested in deriving scaling laws for~$\alpha$ and~$\tau$ as $n \rightarrow \infty$ for given $p_*, p^*, th_*, th^*$, and~$\pi(n)$.  More precisely, we will derive such laws for the median and fixed percentiles of these r.v.s.  This focus on quantiles seems warranted if we want to meaningfully compare our predictions with results of simulation studies such as in~\cite{ahn_terman_2010},  as the mean might be sensitive to outliers such as the ones described in Example~\ref{ex:Landau} of Appendix A that are highly unlikely to ever show up in simulation studies.

Of particular interest is the question whether the parameters $p_*, p^*, th_*, th^*$, and $\pi(n)$ can be chosen in such a way that some form of \emph{chaotic} or \emph{complex} dynamics becomes generic in the corresponding networks.  If we require $\vp  = \vone$, our networks are examples of \emph{Boolean networks,} and it makes sense to adopt the conceptualization of chaotic dynamics proposed by S. A. Kauffman and his followers~\cite{origins}.  In the Boolean context, there are several important hallmarks of chaotic dynamics, two of which are long attractors and  long transients.  In contrast, short attractors and transients are two important hallmarks of
\emph{ordered} dynamics.  Attractors and transients are considered ``long'' in this context if their average lengths, under sampling from a given distribution, increases exponentially, or at least faster than any polynomial, as $n \rightarrow \infty$. As Example~\ref{ex:Landau} shows, the maximum values for~$\alpha$ and~$\tau$ do scale superexponentially, even in the simplest case when $p^* = th^* = 1$, so it is at least conceivable that the same holds  generically for their average value under suitable choices of parameters.  This question appears highly relevant from the point of view of neuroscience, as there may be some randomness to the wiring of neuronal networks in actual organisms, and a number of results in the literature indicate that the potential of  neuronal networks to perform sophisticated computations is highest if their dynamics falls into the \emph{complex regime} that constitutes the boundary between chaos and order  (see, for example, \cite{beggs_timme, friedman_butler,  scarpetta_candia} and references given in these papers).

Preliminary numerical explorations and some results of~\cite{ahn_terman_2010} indicate that for the case~$p^* = th^* = 1$ the most complex dynamics are observed when $\pi(n) \approx \frac{1}{n}$, that is, in the critical window of connection probabilities where the giant component emerges.  Here we will characterize the behavior of~$\alpha$ and~$\tau$ on both sides of the phase transition, and also prove a result that illustrates possible behavior of these systems  inside the critical window. In a planned companion paper  we will explore how these asymptotic theoretic results compare with simulations for bounded~$n$, as well as phase transitions for several related properties of network dynamics.

\section{Statement of results}\label{sec:results}

\medskip

We say that an event~$E_n$ that describes a property of a random structure and depends on parameter~$n$ holds \emph{asymptotically almost surely}
(abbreviated \emph{a.a.s.}), if $P(E_n) \rightarrow 1$ as~$n \rightarrow \infty$ while all other parameters are kept fixed. Similarly,  we say that~$E_n$ occurs
\emph{asymptotically with positive probability} (abbreviated \emph{a.p.p.}),  if $\liminf_{n \rightarrow \infty} P(E_n) > 0$ while all other parameters are kept fixed.

For example, ``$\alpha \geq d$ holds with a.p.p.'' abbreviates the statement that for some fixed probability~$q > 0$ the event~$\alpha \geq d$ occurs with probability~$\geq q$ for all sufficiently large~$n$, while ``$\tau = \theta(\ln(n))$ holds a.a.s.'' abbreviates the statement that for every fixed~$q < 1$ there are constants~$d^l_q, d^u_g$ such that
$P(d^l_q\ln(n) \leq \tau \leq d^u_q\ln(n)) > 1- q$ for all sufficiently large~$n$.

\medskip

Consider connection probabilities of the form $\pi(n) = \frac{c}{n}$, where $c > 0$ is a constant.

For the subcritical case~$c < 1$ we obtained the following characterization.

\begin{theorem}\label{thm:c<1-upper}
Suppose   $0< c <1$  and $\pi(n) = \frac{c}{n}$. Assume $th_* = 1$. Then

\smallskip

\noindent
(a) $\lim_{d \rightarrow \infty} \lim_{n \rightarrow \infty} P(\alpha > d) = 0$.

\smallskip

\noindent
(b)  If $d = 1$ or $d \geq  p_* + 1$, then the equality $\alpha = d$ has a.p.p.
\smallskip

\noindent
(c) The scaling law $\tau = \theta(\ln(n))$ holds a.a.s.
\end{theorem}

By point~(a), the probability of reaching attractors of length~$> d$ dwindles to zero
as~$d$ increases. Thus  the lengths of the attractors scale like a constant, while  the lengths of the transients scale logarithmically.  By~(b), the probability that a periodic attractor will be reached remains strictly between~0 and~1 as $n \rightarrow \infty$.  Proposition~\ref{prop:min-length-periodic} below shows that the theorem gives in fact an exhaustive characterization of possible attractor lengths.

\medskip

For the supercritical case~$c > 1$ we obtained the following two results.

\begin{theorem}\label{thm:c>1-upper-bd}
Suppose  $\pi(n) = \frac{c}{n}$ for some fixed $c > 1$.
 Assume $th_* = 1$.

\smallskip

\noindent
(a)  If $th^* = p^* = 1$, then the inequality $\alpha > 1$ holds a.a.s.

\smallskip

\noindent
(b)  Let  $d \geq p_* + 1$. Then the inequality $\alpha \geq d$ holds a.p.p.

\smallskip

\noindent
(c)  The scaling law $\tau = \Omega(\ln (n))$ holds a.a.s.
\end{theorem}

Thus for~$c > 1$, periodic attractors will be reached a.a.s., at least when $p^* = th^* = 1$.  It follows already from the results of~\cite{JAT} that above  some critical value~$c_{crit}$ the lengths of both attractors and transients scale like a constant.  We conjecture that in the case of attractors, the critical value is~$c_{crit} = 1$, but have not yet been able to prove this.  The following theorem gives a partial result that improves the bounds on~$c_{crit}$ that can be obtained from the approach of~\cite{JAT}.

\begin{theorem}\label{thm:c>ccrit-upper-bd}
 Assume $p_* = p^* = p$, $th^* = 1$ and  $\pi(n) = \frac{c}{n}$ for some fixed $c > p + 1$. Then

\smallskip

\noindent
(a) $\lim_{d \rightarrow \infty} \lim_{n \rightarrow \infty} P(\alpha > d) = 0$.

\smallskip

\noindent
(b)  There exists a constant  $k = k(c,p)$  such that the scaling law $\tau = O(n^{k})$ holds  a.a.s.
\end{theorem}

\medskip

Both in the subcritical and supercritical cases, attractors and transients are relatively short on average, which is indicative of highly ordered dynamics.   Longer attractors  are expected near the phase transition $\pi(n) \sim \frac{1}{n}$.  Our next theorem characterizes $\alpha$ and~$\tau$ for the lower range of the critical window for~$p^* = th^* = 1$.

\begin{theorem}\label{thm:alpha-nonpoly-lower}
Let  $0 < \beta < 1/4$ and let $\pi(n) = \frac{1 - n^{-\beta}}{n}$. Let $1 > \eps > 0$.

\smallskip

\noindent
(a) If $p^* = th^* = 1$, then for every~$k > 0$ the scaling law $\alpha = \omega(n^k)$ holds a.a.s.

\smallskip

\noindent
(b) A.a.s., $\tau = O((\ln n) n^{\beta})$.

\smallskip

\noindent
(c)  If $th^*=1$, then   $\tau = \Omega(n^{\beta})$ holds a.a.s.

\end{theorem}

In particular, Theorem~\ref{thm:alpha-nonpoly-lower} implies that there exists a range of~$\pi(n)$ for which the lengths of attractors scale faster than any polynomial, while the lengths of transients scale like a polynomial.  Thus in this range we observe some, but not all hallmarks of chaotic dynamics, which may indicate that these systems are perched  right on the boundary between order and chaos and exhibit genuinely complex dynamics.  It may seem rather peculiar that for any fixed~$c < 1$ and $\pi(n) = \frac{c}{n}$ the percentiles for~$\alpha$ scale like constants, while the percentiles for~$\tau$ scale like $\log (n)$, while in the window covered by Theorem~\ref{thm:alpha-nonpoly-lower} the median value of $\alpha$ increases much faster than~$\tau$. There is no contradiction here; we can only deduce that the constants implied by Theorem~\ref{thm:c<1-upper} rapidly increase without bound as $c \rightarrow 1^-$.
The  median value of~$\alpha$ grows very slowly though, and the predicted complex behavior is unlikely to show up in simulations for realistic values of~$n$.

\medskip

All results that we have described so far assume~$th_* = 1$ or even~$th^* = 1$.
Networks with $th_* > 1$ behave differently.   The following theorem gives scaling laws for~$\tau$ when $th_* > 1$.
It suggests a different location of the corresponding phase transitions for this case, possibly at $c_{crit} = (th_*!)^{1/th_*}$.

\begin{theorem}\label{thm:th*>1-trans}
Suppose  $\pi(n) = \frac{c}{n}$ for some constant~$c > 0$. Assume
$th_*  > 1$.

\smallskip

\noindent
(a) Then the scaling law $\tau = \Omega(\ln(\ln(n)))$ holds  a.a.s.

\smallskip

\noindent
(b) If $c < (th_*!)^{\frac{1}{th_*}}$, then the scaling law $\tau =  O(\ln(\ln(n)))$ holds a.a.s.
\end{theorem}

The remainder of the paper is organized as follows: In Section~\ref{sec:digraphs} we review properties of Erd\H{o}s-R\'{e}nyi random directed graphs  that are needed for our results on random networks.   Proofs are given in Appendix~B. In Section~\ref{sec:network-tools} we derive properties of network dynamics that are needed for the proofs of  Theorems~\ref{thm:c<1-upper}--\ref{thm:alpha-nonpoly-lower}, which are presented in  Section~\ref{sec:proofs}. The proof of Theorem~\ref{thm:th*>1-trans} will be given in Appendix~C as it uses slightly different techniques. Finally, Section~\ref{sec:discussion} lists some open problems for further investigation.

\section{Erd\H{o}s-R\'{e}nyi random digraphs}\label{sec:digraphs}

Let $D$ be a digraph on~$[n]$, and let $b >0$ be a constant.  Following~\cite{Karp} we say that a set of nodes~$V \subseteq [n]$ is \emph{$b$-small} if
$|V| < b \ln n$ and that $V$ is \emph{$b$-large} otherwise. We call a digraph~$D$ \emph{supersimple}  if it
contains at most one directed cycle, and, moreover, if~$C$ is a directed cycle in~$D$ and~$j$ is a node outside of~$C$, there are no two distinct  directed paths~$Pt_1, Pt_2$ from $C$ to~$j$ such that both~$Pt_1$ and~$Pt_2$ are arc-disjoint from~$C$ and there are no two  directed paths~$Pt_1, Pt_2$ from~$j$ to~$C$  such that both~$Pt_1$ and~$Pt_2$ are arc-disjoint from~$C$.

 For a set~$V \subseteq [n]$ let $CYC_{[\ell]}(V)$ denote the property that~$D$ contains a directed cycle of length~$\ell$ all of whose nodes are inside~$V$. Instead of $CYC_{[\ell]}([n])$ we write $CYC_{[\ell]}$.  We let $CYC$ denote the property that~$D$ has at least one directed cycle.

 For  positive integers~$I, J$ and a positive constant $\kappa < 1$, let $P2C(I, J, \kappa)$ denote the property that there exists a set
$P2P = \{q_r\ell_r: \, r \in [J]\}$ such that $q_r < I$ is an odd prime,  $\ell_r > n^\kappa$ is prime, the numbers $\ell_r$ are pairwise distinct, and for all $r \in [J]$ there exists a directed cycle of length~$q_r\ell_r$ in~$D$.

Let $CYCU$ denote the set of all nodes~$i$ such that~$UC(i)$ contains a directed cycle.

We let $\varrho(c)$ denote the  unique solution of the equation

\begin{equation}\label{eqn:first-entry-proportion-giant-size}
e^{-c\varrho} = 1 - \varrho
\end{equation}
in the interval~$(0,1)$.

The following lemma summarizes the properties of Erd\H{o}s-R\'{e}nyi random digraphs that we will need for deriving the main results of this paper.  For the proof, see
Appendix~B.

\begin{lemma}\label{lem:Karp}
Let~$D$ be a digraph on~$[n]$ whose arcs are randomly and independently drawn with probability $\pi(n)$.

\begin{itemize}
\item[(A)] Assume $\pi(n) = \frac{c}{n}$ for some constant $0 < c < 1$. Then  there exist  positive constants $a = a(c)$ and $b = b(c)$ such that

\begin{itemize}
\item[(A1)] A.a.s., all sets $UC(i), DC(i)$ will be $b$-small and supersimple.
\item[(A2)] A.a.s, there exists $i \notin CYCU$ such that $L_{\max}^s(i) \geq a \ln n$.
\item[(A3)] $\lim_{k \rightarrow \infty} \lim_{n \rightarrow \infty} P(|CYCU| \geq k) = 0$.
\item[(A4)] For each $\ell > 1$, property $CYC_{[\ell]}$ has a.p.p. and $\neg CYC$ also has a.p.p.
\end{itemize}

\item[(B)] Let $\pi(n) = \frac{1 - n^{-\beta}}{n}$ for some constant  $0 < \beta < 1/4$ and let~$\eps > 0$. Then

\begin{itemize}
\item[(B1)] A.a.s., all sets $UC(i), DC(i)$ will be supersimple.
\item[(B2)] $\forall J, \kappa \, \exists I \     \lim_{n \rightarrow \infty} P(P2C(I, J, \kappa))  > 1 - \eps$.
\item[(B3)] A.a.s. $L_{\max} \leq (1 + \eps)(\ln n) n^{\beta}$ and $L_{\max}^s \geq (1 - \eps)n^{\beta}$.
\item[(B4)] $\lim_{t \rightarrow \infty} \lim_{n \rightarrow \infty} P(L_{\max}(1) \geq t) = 0$.
\item[(B5)] For all fixed~$t$, with a.p.p. we have $L_{\max}^s(1) \geq t$.
\end{itemize}

\item[(C)] If $\pi(n) = \frac{c}{n}$ for some constant $c > 1$, then there exist positive constants $a = a(c)$ and $b = b(c)$ such that

\begin{itemize}
\item[(C1)] A.a.s., $|GC| = (\varrho(c)^2 + o(1))n$ and $|DG| = (\varrho(c) + o(1))n$.
\item[(C2)] A.a.s., for all $i \notin DG$, the set  $UC(i)$ is $b$-small and supersimple.
\item[(C3)]  $\lim_{k \rightarrow \infty} \lim_{n \rightarrow \infty} P(|CYCU \backslash DG| \geq k) = 0$.
\item[(C4)] A.a.s.,  there exists $i \notin CYCU$ such that  $L_{\max}^s(i) \geq a \ln n$.
\item[(C5)] A.a.s., for every $i \in UG$  and $j \in DG$ there exists a directed path from~$i$ to~$j$ of
length $\leq b\ln n$.
\end{itemize}
\end{itemize}
\end{lemma}

\section{Properties of network dynamics}\label{sec:network-tools}

In this section we prove results about network dynamics that are basic tools for deriving our main theorems. In order to reduce clutter, we will usually omit the standing assumption that a network~$N ~= ~ < D, \vec{p}, \vec{th}>$ on~$[n]$ and an initial state~$\vs(0)$ are given. Let us begin with an elementary observation

\begin{proposition}\label{prop:min-length-periodic}
Any periodic attractor must have length $\alpha \geq p_* + 1$.
\end{proposition}

\noindent
\textbf{Proof:} In any state of a periodic attractor, some node~$i$ must fire.  But if $s_i(t) = 0$, then we have $s_i(t+1) = 1, \dots , s_i(t+p_i) = p_i$.  Since $p_i \geq p_*$, it follows that a periodic attractor must contain at least $p_*+1$ distinct states. $\Box$

\subsection{Decomposition of network dynamics}\label{subsec:decomposition}

Let us make another simple but crucial observation: If $j \in UC(i)$, then $s_j(t+1)$ depends only on the states $s_k(t)$ such that $k = j$ or $<k, j> \, \in A_D$.  In particular, such that $k \in UC(i)$. By induction, this implies the following.

\begin{proposition}\label{prop:upstream}
Let  $\vs(0), \vsst(0)$ be two initial conditions, and let $i \in [n]$.
If $s_j(0) = s^*_j(0)$ for all $j \in UC(i)$, then $s_j(t) = s^*_j(t)$ for all $t \geq 0$ and $j \in UC(i)$, in particular, $s_i(t) = s^*_i(t)$ for all $t \geq 0$.
\end{proposition}

We will refer to the trajectories of the nodes in~$UC(i)$ as the \emph{internal dynamics} of $UC(i)$, and let $\alpha_i, \tau_i$ denote the lengths of the attractor and the transient in the internal dynamics of~$UC(i)$. Then Proposition~\ref{prop:upstream} implies

\begin{equation}\label{eqn:alpha-tau-upstream}
\alpha = lcm\{\alpha_i: \ i \in [n]\} \qquad \mbox{and} \qquad \tau = \max\{\tau_i: \ i \in [n]\}.
\end{equation}

\subsection{Necessary and sufficient conditions for long transients}\label{subsec:transients}

Example~\ref{ex:Landau} of Appendix~A shows that long transients may result from the interactions of nodes in many relatively short cycles. In contrast, if~$D$ contains at most one directed cycle, transients cannot be much longer than the longest directed path. The following two propositions make this connection explicit and also give some information about possible attractor lengths when~$D$ is simple.

\begin{proposition}\label{prop:too-short-cycles}
 Suppose $UC(i)$ is acyclic.  Then $\alpha_i = 1$ and
 $\tau_i \leq L_{\max}(i) + p^* \leq |UC(i)| + p^* - 1$.
\end{proposition}

\medskip
\noindent \textbf{Proof:}  Consider a node $j \in UC(i)$ and assume $s_j(t^*) = 0$ for some~$t^* > 0$. Then there exists a sequence $( j_{t^*}, \ldots ,j_0)$ of nodes in~$UC(i)$ with $j_{t^*} = j$, $s_{j_t}(t) = 0$, and $<j_{t}, j_{t+1}> \, \in A_D$ for all $t < t^*$.  This sequence is in general not unique; we will call each such sequence \emph{a history of $s_j(t^*) = 0$.}  If
$UC(i)$ is acyclic, then each history corresponds to a directed path of length~$t^*$ in~$UC(i)$, and consists of $t^* + 1$ pairwise distinct nodes.  It follows that  we must have
$s_j(L_{\max}(i) + 1) \geq 1$ for all $j \in UC(i)$, and the steady-state attractor will be reached after at most~$p^* -1$ additional steps. Thus  $\alpha_i = 1$ and $\tau_i \leq L_{\max}(i) + p^*$.  Since the number of arcs in any relevant directed path cannot exceed $|UC(i)|$, the  inequality $L_{\max}(i) + p^* \leq |UC(i)| + p^* - 1$ always holds.
$\Box$

\begin{proposition}\label{prop:only-1-cycle-full}
 Suppose $D$ contains exactly one directed cycle~$C$.

\smallskip

\noindent
(a) If $\alpha > 1$, then $gcd \{\alpha, |C|\} > 1$.

\smallskip

\noindent
(b) If $D$ is supersimple, then  $\tau \leq  2 L_{\max} + ((p^*+1)^{p^*} + 1)|C| + 4 p^* - 3  = O(L_{\max})$.
\end{proposition}

\medskip

\noindent
\textbf{Proof:}  Let $V_u = [n] \backslash DC(C)$, and let  $L_{\max}^u$ be the length of the longest directed path that consists of nodes in $V_u$. Then the restriction of~$D$ to nodes in~$V_u$ is acyclic, and it follows from Proposition~\ref{prop:too-short-cycles} that at time  $t_0 =  L_{\max}^u  + p^*$ all nodes in $UC(i)$ for $i \in V_u$ will have reached the end of their refractory period and
will not subsequently fire. In particular, they will not influence the internal dynamics of $UC(j)$ at any time  $t \geq t_0$ for any of the nodes~$j \in DC(C)$.  Thus for the remainder of this argument we may disregard  these nodes and consider what happens in the system restricted to $DC(C)$ for the initial state~$\vs(t_0)$.

In the restricted system we will have $UC(i) = C$ for each~$i \in C$, and  Proposition~1 and Theorem~3 of~\cite{AJ}  imply  that the length of the attractor  is a divisor of~$|C|$ and the length of the transient is bounded from above by $|C| + 2p^* - 3$. Thus at time $t_1 = t_0 + |C| + 2p^* - 3$ our original system will have reached a state in the attractor of the internal dynamics of~$UC(i)$ for $i \in C$.  Then for some $\alpha_C$ we will have $\alpha_i = \alpha_C$ for every~$i \in C$, and  the number $\alpha_C$ is a divisor of~$|C|$.  Thus $\alpha_C$ divides both~$|C|$  and~$\alpha$
by~(\ref{eqn:alpha-tau-upstream}). If $\alpha_C > 1$, this implies~(a). If $\alpha_C = 1$, then we will also have~$\alpha = 1$ (see Case~1 below), and point~(a) vacuously holds.

\medskip

For point~(b),  first note that the order-of-magnitude estimate follows from the fact that~$|C| \leq L_{\max}$.  For the proof of the first inequality, we distinguish two cases.

\medskip

\noindent
\emph{Case~1: $\alpha_C = 1$}

\smallskip

In this case $s_j(t_1) = p_j$ for all nodes in $V_u \cup C$.  Since the restriction of~$D$ to the nodes in $DC(i) \backslash C$ is acyclic, Proposition~\ref{prop:too-short-cycles} applies, and $\vs(t_1 + L_{\max}^d + p^*) = \vp$, where $L_{\max}^d$ is the length of the longest directed path in $[n] \backslash (V_u \cup C)$.  It follows that

\begin{equation}\label{eqn:gamma-u-C-gamma-d}
\tau \leq  L_{\max}^u + 4p^*  + |C|  - 3 + L_{\max}^d \leq  2 L_{\max} + |C| + 4p^* - 3.
\end{equation}

\medskip

\noindent
\emph{Case~2: $\alpha_C > 1$}

\smallskip

For  $d > 0$, let $L_d$ be the set of nodes $j \in DC(C)$ such that there exists a directed path~$Pt = (i = i_0, i_1, \ldots i_{d-1}, j)$ of length~$d$ from some $i \in C$ to $j$ with  $\{i_1, \ldots, i_{d-1}, j\} \subset DC(C) \backslash C$. By assumption, this directed path is unique, so node $j \in L_d$ can take input only from exactly one~$i_{d-1} \in L_{d-1}$  (or from $i \in C$ if $d =1$) and possibly some nodes in~$V_u$ that  never fire at any time $t \geq t_0$ and can be ignored.

First consider a node~$j \in L_1$.   Then there are  at most $(p_j+1)|C|$ distinct states that the nodes in the set $C \cup \{j\}$ can take at times~$t \geq t_1$, and
it follows that~$\tau_j < t_1 + (p_j+1)|C|$ and $\alpha_j \leq (p_j+1)\alpha_{i} \leq (p_j + 1)|C| \leq (p^* + 1)|C|$.

Now consider~$j \in L_d$  for $d > 1$ that takes input from~$i_{d-1}  \in L_{d-1}$. Node~$j$ can fire at time
$t \geq \max \{1, t_1 - p^* -1\}$ only if node~$i_{d-1}$ fires at time~$t-1$.

If $s_{i_{d-1} }(t-1) = 0$ for~$t \geq t_1$ and $p_j \leq p_{ i_{d-1} }$, then we must have   $s_{j} (t-1) = p_j$,  as node~$j$ can fire at any time $\geq t_1 - p^* -1$ only if node~$ i_{d-1}$ fires at the previous time step.  Thus the implication $s_{ i_{d-1} }(t-1) = 0 \Rightarrow s_j(t)  = 0$ will hold for all $t \geq t_1$, and it follows that~$\tau_j \leq \max\{t_1, \tau_{ i_{d-1} } + 1\}$ and $\alpha_j = \alpha_{ i_{d-1} }$.

If $s_{ i_{d-1} } (t-1) = 0$ for~$t \geq t_1$ and $p_j > p_{ i_{d-1} }$, then we need to take into account the possibility that  node~$j$ will not yet have reached the end of its refractory period at time~$t-1$. Think about this as the worst-case scenario.
However, a simple inductive argument shows that for $t \geq d + t_1$ the latter can only occur if there is no node~$i_k \in Pt$ with $k < d-1$ such that $p_{i_k} \geq p_j$,   as otherwise node~$i$ could not have fired twice in the interval $[t - 1 - p_j, t-1]$. It follows that the worst-case scenario can occur at most $p^* -1$ times along any directed path~$Pt$ as above.  If the worst-case scenario does occur, then the same state of $UC(j)$ must occur twice in the interval
$[t_1, t_1 + (p_j+1)\alpha_{ i_{d-1} }]$.

Now it follows by induction over~$d$ that $\tau \leq t_1 + d + (p^* + 1)^{p^*}|C|$,  which implies~(b).
 $\Box$.

\begin{remark}
We cannot in general conclude that $\alpha$ is a divisor of~$|C|$.  To see this, consider $D$ on~$[3]$ with arc set $A_D = \{<1,2>, <2,1>, <1,3>\}$, $\vth = \vone$, $p_1 = p_2 = 1$, $p_3 = 2$, $s_1 = 0$,  $s_2=s_3=1$.  Then the attractor has length~$4$, whereas the only directed cycle has length~$2$. We do not know whether the assumption that~$D$ is supersimple is needed in part~(b), but we will apply the result only in situations where this property holds.
\end{remark}

The existence of long directed paths does not all by itself imply the existence of long transients; consider for example an initial state that is a constant vector.
However, under some additional assumptions  the existence of long directed paths  guarantees that long transients will occur with high probability.

\begin{lemma}\label{lem:path-implies-transient}
Assume $\frac{c^-}{n} \leq \pi(n) \leq \frac{c^+}{n}$, and let $Pt = (i_0, i_1, \ldots, i_\ell)$ be a sequence of pairwise distinct nodes with $th_{i_j} = 1$ for all
$j\in \{0, \ldots , \ell\}$.  Then for all $\delta > 0$ there exist  a constant  $R = R(p_*, p^*, c^-, c^+, \delta) > 0$ that does not depend on~$n$ or~$\ell$  such that
\begin{equation}\label{eqn:P-cond-path-to-trans}
P(\tau \geq R \ell \, | \, D \ \mbox{is  supersimple} \, \& \, Pt \  \mbox{is a straight directed path} \, \& \, |UC(i_0)| < \frac{n}{2}) > 1 - \delta.
\end{equation}
\end{lemma}

\noindent
\textbf{Proof:} Let us assume that~$Pt = (i_0, \ldots , i_\ell)$ is a straight directed path in~$D$. If the subdigraph of~$D$ that is induced by $UC(i_\ell)$ is  not acyclic, then there exists a largest node~$i_k$ that is downstream from any directed cycle~$C$.  If~$D$ is supersimple, then both~$C$ and the path from~$C$ to~$Pt$ are unique.   This makes~$k$ a well-defined r.v.  Since~$Pt$ is assumed straight, it must be disjoint from~$C$, so that we have a symmetric situation which guarantees that~$k$ is uniformly distributed over the interval~$[0, \ell]$.
It follows that for $k = \lfloor \frac{\delta \ell}{2}\rfloor$ with probability $\geq 1 - \frac{\delta}{2}$ the subdigraph of~$D$ that is induced by $UC(i_{k-1})$ is  acyclic.

Now assume the latter and consider  the path
$Pt_1 = (i_0, \ldots , i_{k-1})$  instead of~$Pt$.  Node~$i_{k-1}$  will eventually stop firing. Thus the time when the last firing of node~$i_{k-1}$ occurs must be $\leq \tau$
and we only need to derive  a lower bound on the expected time of this last firing.

Let us first derive an estimate that works for the special case $p^* = 1$ and illustrates important ideas. Later we will give an argument that works for the general case.

Let $t_j$ be the  last time~$t$ when nodes $i_{j}$ fires.  If node~$i_j$ never fires, we define~$t_j = -1$.
If node~$i_j$  fires at time~$t_j$, then node~$i_{j+1}$  will fire at time   $t_j+1$   unless it also fires at time~$t_j$. It follows that for all~$j<k-1$  we have

\begin{equation}\label{eqn:tj-ineqs-basic}
t_{j+1} \geq t_j \geq -1.
\end{equation}

Moreover, if node~$i_{j+1}$ has indegree~$1$, and~$t_j > 0$, then it cannot fire at time~$t_j$, as the only node capable of inducing such firing would be~$i_j$.  Thus we must have $t_{j+1} = t_j+ 1$ in this case. For large enough~$u$, the probability that $t_j > 0$ for at least one of the nodes $i_0, \dots , i_u$ is arbitrarily close to~1. It follows by induction that if $DO$ denotes the number of indices~$j$ with~$u \leq j \leq k-1$  such that node~$i_j$ has indegree~1, then  $t_{k} \geq DO$. For each relevant~$j$ the conditional probability of the event that~$i_j$ has indegree~1 given that~$Pt$ has the properties specified above is larger than a positive constant~$q = q(c)$ that does not depend on~$n$. This implies the lemma for any choice of~$R$ with $R < q\lfloor \frac{\delta}{2}\rfloor$.

Now let us prove the result in the general case. Consider any directed path~$Pt^* = (j_0, \ldots , j_u)$ in~$D$.  We say that \emph{$Pt^*$ has a forcing extension} if

\begin{itemize}
\item $s_{j_0}(0) = 0$,
\item for $w \in [2p^*]$ node~$j_w$ has indegree~1,
\item for $w \in [2p^*]$ the equality $s_{j_w}(0) = p_{j_w}$ holds,  and
\item for $w > 2p^*$ there does not exist a directed path that ends at~$j_w$, has length  $> w - 2p^*$   and does not contain node~$j_{2p^*}$.
\end{itemize}

The key observation here is that if $Pt^*$ has a forcing extension, then for all $w \in [u]$, node~$j_w$ will fire at time~$t = w$.  For   $w \leq 2p^*$  this follows immediately from the choice of the initial condition; for  $w > 2p^*$  this can easily be shown by induction over~$w$:
Node $j_w$ will receive its last firing input from \emph{outside of}~$Pt^*$ at time $t < w - 2p^*$, and the specifications of the initial condition implies that it will not receive any firing input during the time interval $[w- 2p^*, w -1)$.  Thus  $s_{j_w}(w-1) = p_{j_w}$,  and the firing of node~$j_{w-1}$ at time~$w-1$ will induce a firing of node~$j_w$ at time~$w$.

Now assume $j$ is any node such that~$UC(j)$ is acyclic. Let $r$ denote the  probability that the longest directed path that ends at~$j$ has a forcing extension. We claim that~$r$ is bounded from below by a constant~$W = W(c^-, c^+, p^*) > 0$ that does not depend on~$n$.  In order to see this, assume that $Pt^-$ is a straight path of that ends at~$j$ and has the property that every directed path that ends at~$j$ and has length larger than the length of~$Pt^-$ contains~$Pt^-$. Assume moreover
that~$|UC(j)| < \frac{n}{2}$. Then the probability that
$Pt^-$ can be extended to a straight path with a forcing extension, conditioned on all the above assumptions on~$Pt^-$, is clearly positive and bounded from below by a constant~$W$ that does not depend on~$n$.

Now let~$Pt$ be as at the beginning of the proof of the lemma.  Fix an integer~$Sg$ and partition~$Pt$ into~$Sg$ consecutive segments of length $\geq \lfloor\frac{\ell - Sg}{Sg}\rfloor$ each. Assume~$Pt_{xy} = (i_x, \ldots , i_y)$ is a given such segment.  Then the probability that either~$Pt_{xy}$ itself has a forcing extension or there
exists a directed path~$Pt^*_{xy}$ of length  $\geq y - x + 2p^*$  with a forcing extension that branches off~$Pt_{xy}$ at some point~$i_z$ is also $\geq r$.
By symmetry, the probability that $p_{i_z} = p^* = s_{i_z}(t)$, where $t$ is equal to the length of~$Pt^*$ minus~$y - z$ is at least~$(p^*+1)^{-2}$. Thus with
probability~$\geq (p^*+1)^{-2}r$ node $i_y$ will fire at some time $t \geq \lfloor\frac{\ell - Sg}{Sg}\rfloor - 2p^*$.  If we consider only the set of even-numbered segments, the events that the terminal nodes fire at times as specified will be independent, and the result follows by choosing~$Sg$ sufficiently large relative to~$\alpha$.
$\Box$

\subsection{Sufficient conditions for long attractors}\label{subsec:cond-attractors}

Proposition~\ref{prop:too-short-cycles} shows that the existence of at least one directed cycle in~$D$ is a necessary condition for the existence of attractors of length~$>1$.  Here we derive some sufficient conditions for the existence of attractors of length~$\geq d$ for certain positive integers~$d \geq 2$.

Existence of directed cycles in~$D$  all by itself does not guarantee that a periodic attractor will be reached. This can easily be seen by considering the case where~$D$ is a directed cycle and $\vs(0)$ is constant.  It turns out though that when $p^* = th^* = 1$, this is in a sense the only counterexample.
More precisely,  let~$C = (i_1, \ldots , i_k, i_{k+1} = i_1)$ be a directed cycle in~$D$. Define   $S_C = \{\vs \in St_N: \exists i_j, i_\ell \in C \ s_{i_j} = 0 \ \& \ s_{i_\ell} = 1\}$.    The following result is a straightforward generalization of Proposition~6.23 of~\cite{bookchapter}.

\begin{proposition}\label{prop:S_C-is-invariant:online}
Assume $p_i = th_i = 1$ for all $i \in C$. Then the set~$S_C$ is forward invariant under the dynamics of~$N$, that is, if $\vs(t) \in S_C$, then $\vs(t+1)\in S_C$.
In particular, $S_C$ is disjoint from the basin of attraction of the steady state attractor.
\end{proposition}

\noindent
\textbf{Proof:} Let $\vs(t) \in S_C$ and let $i_j, i_\ell \in C$ be such that $s_{i_j}(t) = 0$ and $s_{i_\ell}(t) = 1$. Wlog $j < \ell$ and there must be an index $m$ with $j \leq m < \ell$ such that $s_{i_m}(t) = 0$ and $s_{i_{m+1}}(t) = 1$.    Then  $s_{i_m} (t+1) = 1$ and $s_{i_{m+1}} (t+1) = 0$,
and thus $\vs(t+1) \in S_C$.  The second sentence of the proposition follows from the first and the fact that the steady state $\vec{p}$ is not in~$S_C$. $\Box$

\bigskip

Let $NSC_{[\ell]}$ denote the property that there exists a directed cycle~$C$ of length~$\ell$, composed exclusively of nodes~$j$ with $p_j = th_j = 1$, such that the initial condition~$\vs(0)$ takes both values~$0$ and~$1$ on~$C$. This condition makes only sense if $p_* = th_* = 1$. We will use witnesses of property~$NSC_{[\ell]}$ in the proof of Theorem~\ref{thm:alpha-nonpoly-lower}.

When $p_* > 1 = th_*$, things become a bit more complicated, as we are no longer automatically assured that a node in the cycle that does not fire is at the end of its refractory period. We do not know of a precise analogue of Proposition~\ref{prop:S_C-is-invariant:online} for this case, but we can at least formulate a condition for the case when there are no additional arcs with targets in the cycle.

Let $NSC_{[\ell]}^1$ denote the property that there exists a directed cycle~$C$ of length~$\ell$, all of whose nodes~$j$ have  indegree~1,   refractory period~$p_j = p_*$, and firing threshold $th_j = 1$, such that the initial condition~$\vs(0)$ satisfies $s_i(0) = 0$ for exactly one node~$i \in C$ and
$s_j(0) = p_j$ for all other nodes~$j \in C$.

 In a cycle~$C$ that witnesses~$NSC_{[\ell]}^1$, at any given time~$t$ at most one node can fire.  Moreover, if the cycle is sufficiently long, after having fired, a node will have reached the end of its refractory period when the firings have traveled around the circle, and will fire again. Thus~$\alpha_i = \ell$ for all~$i \in |C|$. This observation, together with
Proposition~\ref{prop:S_C-is-invariant:online}, implies that

\begin{equation}\label{eqn:NSC-implies-alpha>1}
\begin{split}
&(\exists \ell \geq 2 \ NSC_{[\ell]}) \Rightarrow \alpha > 1,\\
& NSC_{[\ell]}^1 \Rightarrow \alpha \ \mbox{is a multiple of} \   \ell.
\end{split}
\end{equation}

Some special initial states allow us to drop the requirement that all nodes in~$C$ have indegree~$1$.
Let $NSC_{[\ell]}(p)$ denote the property that there exists a directed cycle~$C = (i_0, i_1, \dots , i_\ell = i_0)$,  composed exclusively of nodes~$j$ with $p_j = p$ and $th_j = 1$, such that $s_{i_j}(0) = j \ mod \ (p+1)$ for all $i_j \in C$. This property implies that $\ell$ is a multiple of~$p+1$.  If a cycle~$C$ witnesses~$NSC_{[\ell]}(p)$, then all nodes will receive firing inputs from \emph{within} the cycle exactly at times when they reach the end of their refractory periods, and input from outside the cycle becomes irrelevant.  This implies the following analogue of~\eqref{eqn:NSC-implies-alpha>1}:

\begin{equation}\label{eqn:NSCell(p)-implies-alpha>1}
(\exists \ell, p  \ NSC_{[\ell]}(p)) \Rightarrow   \alpha \ \mbox{is a multiple of} \ p+1.
\end{equation}

\begin{lemma}\label{lem:NSC-c-fixed}
 Assume $th_* = 1$ and  $\pi(n) \geq \frac{c}{n}$ for some fixed $c > 0$. Then

\smallskip

\noindent
(a) For every $\ell \geq  p_*+1$ property $NSC_{[\ell]}^1$ holds  a.p.p.

\smallskip

\noindent
(b) If $p_* \leq p \leq  p^*$ and  $c > (p+1)(p^* - p_* + 1)th^*$, then for every integer~$k$  a.a.s there exists $\ell \geq k$ such that property
 $NSC_{[\ell]}(p)$  holds.
\end{lemma}

\noindent
\textbf{Proof:} Assume~$\ell \geq p_* + 1$. Then  $CYC_{[\ell]}$  has a.p.p. For $c < 1$ this  follow directly from Lemma~\ref{lem:Karp}(A4); for $c \geq 1$ this follows from the observation that   $P(CYC_{[\ell]})$ increases monotonically with respect to~$c$. Consider a sequence $C = (i_0, i_1, \dots , i_\ell = i_0)$. The conditional probability that~$C$ witnesses~$NSC_{[\ell]}^1$, \emph{given} that all the arcs that make~$C$ a directed cycle are included in~$A_D$, is approximately equal to   $\frac{\ell e^{-c\, \ell}}{(p^* - p_* +1)^\ell (th^*)^\ell}$, as there are $\ell$ possible locations of the node that fires initially, the refractory period and initial state of each node in~$C$ are determined by condition $NSC_{[\ell]}^1$ and the location of the initial firing, and the number of all arcs with targets in~$C$, \emph{other} than the arcs of the directed cycle, has approximately a Poisson distribution with parameter~$c\, \ell$. Thus the conditional probability   is bounded from below by a fixed positive constant that does not depend on~$n$.  This implies point~(a).

\medskip

 Part~(b) is phrased here in more general form than we strictly need for our results as the proof is almost identical to the one for the special case  $p^* = p_* = p$ for which we will use the lemma.  It will be given in Appendix~B as it depends on some parts of the proof of Lemma~\ref{lem:Karp}. $\Box$

\subsection{Minimally cycling nodes and eventually minimally cycling nodes}\label{subsec:mincyc}

\begin{definition}\label{def:EV-min-cyc}
 Consider $t_0 < t_1 \leq \infty$. We call $t \in [t_0, t_1)$ \emph{an interval of uninterrupted firings of node~$i$} if there are no times $t$ with $t_0 \leq t < t_1$ such that  $s_i(t) = s_i(t+1) = p_i$.   A node~$i$ is \emph{eventually minimally cycling} if there exists a $t_0 \geq 0$ such that $[t_0, \infty)$ is an interval of uninterrupted firings of node~$i$. The smallest time~$t_0$ for which this happens will be referred to as \emph{the time when node~$i$ becomes minimally cycling.} We say that~$i$ \emph{is minimally cycling} if $t_0 = 0$.
\end{definition}

These notions behave nicely with respect to downstream components of~$D$.

\begin{lemma}\label{lem:ev-min-cyc}
Assume  $p_* = p^* = p$ and $th^* = 1$ and suppose $j$ is a minimally cycling node.

\smallskip

\noindent
(a) Let  $i \in DC(j)$. Then $i$ also is eventually minimally cycling.

\smallskip

\noindent
(b) Let~$L$ be such that for every $i \in DC(j)$ there exists a directed path of length $\leq L$ from~$j$ to~$i$.
Moreover, let $UC^* = UC(DC(j)) \backslash DC(j)$, let $\tau^* = \max \{\tau_k: \ k \in UC^*\}$,
$\alpha^* = lcm(\{\alpha_k: \ k \in UC^*\})$.  Then the smallest time $\tau^{all}$ at which every  node~$i \in DC(j)$ becomes minimally cycling satisfies

\begin{equation}\label{eqn:tmc(j)}
\tau^{all} \leq \tau^* + L + \alpha^*(p+1)^{L+1}.
\end{equation}
\end{lemma}

\noindent
\textbf{Proof:} Part~(a) follows from the last sentence of the following proposition by induction on the length of the shortest path from~$j$ to~$i$.

\begin{proposition}\label{prop:ev-min-cyc2}
Assume $th_i = 1$ and suppose $i, j$ are such that $<j, i> \, \in A_D$ and $p_i+1$ is divisible by $p_j+1$. If $[t_0, t_1)$ is an interval of uninterrupted firings of node~$j$, then there can be at most
$p_j$ times $t \in [t_0 + 1, t_1+1)$ such that $s_i(t) = s_i(t+1) = p_i$.  In particular, if node~$j$ is eventually minimally cycling, then so is node~$i$.
\end{proposition}

\smallskip

\noindent
\textbf{Proof:} Let $i, j, [t_0, t_1)$ be as in the assumptions. For simplicity of notation, assume that $s_j$ fires at all times $t \in [t_0, t_1)$ such that $t$ is divisible by $p_j+1$. Consider the  times $t_\theta \in [t_0+1, t_1+1)$ with $s_i(t_\theta)=s_i(t_\theta+1)=p_i.$  Define $r(\theta)$ as the remainder of $t_\theta$ under division by $p_j+1.$ Note that for each $\theta$ we must have $t_{\theta+1} = t_\theta +1+q(p_i+1)$ for some nonnegative integer $q$. Since $p_j+1$ divides $p_i+1$, it follows that   $r(\theta+1) =  r(\theta)+1$   for all $\theta$, and if there were more than~$p_j$ distinct times~$t_\theta$, at least one of them would satisfy   $t_\theta = 0 \ mod~  (p_j+1)$. Then $r_\theta=0$. This  means that at time $t_\theta$ we would have $s_j(t_\theta)=0$ and $s_i(t_\theta)=p_i$. Thus $s_i(t_\theta+1)=0 \neq p_i$, which contradicts the definition of $t_\theta$.
 $\Box$

\bigskip

For the proof of part~(b), let us first make the following observation:

\begin{proposition}\label{prop:AUTO-stretch}
Let the notation be as in Lemma~\ref{lem:ev-min-cyc}, and suppose that $t_0, t_1$ are times such that
$t_0 > \tau^*, t_1 > t_0 + \alpha^*(p + 1)$, and  $[t_0, t_1)$ is an interval of uninterrupted firing for every node~$i \in DC(j)$. Then every node~$i \in DC(j)$ becomes  already minimally cycling no later than at time~$t_0$.
\end{proposition}

\noindent
\textbf{Proof:} The assumptions imply that the state of the system restricted to $UC(DC(i))$ at time~$t_0 +  \alpha^*(p + 1)$ is the same as at time~$t_0$, which means that the restricted system has already reached its attractor at time~$t_0$. All nodes in~$DC(j)$ are minimally cycling in the attractor, and the result follows. $\Box$

\bigskip

Now let~$j$ be a minimally cycling  node as in the assumption of the lemma. Then   $[t_0, t_1) = [\tau^*, \tau^* + \alpha^*(p+1)^{L+1})$ is an interval of uninterrupted firings of node~$j$. Proposition~\ref{prop:ev-min-cyc2} shows that if  $<j, i> \, \in A_D$, then $[t_0, t_1)$ is a union of at most~$p+1$ subintervals of uninterrupted firings of node~$i$, and at least one of them must have length $\geq \frac{t_1 - t_0}{p+1}$. It follows by induction over the length of the shortest path from~$j$ to~$i$ that every node~$i \in DC(j)$ has an  interval of uninterrupted firings $[t_0^i, t_1^i)$ of length $\geq \alpha^*(p + 1)$ with $t_0^i \leq \tau^* + L + \alpha^*(p+1)^{L+1}$, and part~(b) follows from Proposition~\ref{prop:AUTO-stretch}. $\Box$

\bigskip

Let
$MNODE(g)$ (respectively: $EVMNODE(g)$) denote the property that there exists a giant strongly connected component~$GC$ in~$D$ and the trajectory of~$\vs(0)$ contains at least one (eventually) minimally cycling node in~$GC$.  The next lemma is given here in more general form than we strictly need for our results as the proof is almost identical to the one for the special case $p^* = p_* = p$ and $th^* = 1$ for which we will use the result.

\begin{lemma}\label{lem:P(evmnode(g))}
Suppose  $th_* = 1$  and $\pi(n) = \frac{c}{n}$ for some fixed   $c > (p_*+1)(p^* - p_* + 1)th^*$. Then $MNODE(g)$ holds a.a.s.
\end{lemma}

\noindent
\textbf{Proof:}  Let $q$ be any fixed positive probability. Choose~$k$  sufficiently large so that for all sufficiently large~$n$ the  probability that there are at least~$k$ vertices outside of $DG$ that belong to a directed cycle in~$D$ is less than~$\frac{q}{3}$.    Lemma~\ref{lem:Karp}(C3) implies that such~$k$ exists. By symmetry,  for all sufficiently large~$n$ the  probability that there are at least~$k$ vertices outside of $UG$ that belong to a directed cycle in~$D$ is also less than~$\frac{q}{3}$.

Now Lemma~\ref{lem:NSC-c-fixed}(b) implies that for sufficiently large~$n$, with probability $> 1 - \frac{q}{3}$ there will exist some~$\ell > k$ such that property~$NSC_{[\ell]}(p_*)$ holds.  But if
$C$ is a directed cycle that witnesses property~$NSC_{[\ell]}(p_*)$, then every node in~$C$ is minimally cycling.  Moreover, $C$ must contain at least~$k$ nodes, and with probability $> 1 - \frac{2q}{3}$ at least one (and hence every) node in~$C$ will be in $GC$.  This implies
 $P(MNODE(g)) > 1-q$ and proves the lemma.  $\Box$

\section{Proofs of Theorems~\ref{thm:c<1-upper}--\ref{thm:alpha-nonpoly-lower}}\label{sec:proofs}

\subsection{Proof of Theorem~\ref{thm:c<1-upper}}\label{subsec:PfThm1}

\medskip

Consider $i \in [n]$. If $UC(i)$ is acyclic, then~$\alpha_i = 1$.  Otherwise
$\alpha_i$ is trivially bounded by the size of the state space of the internal dynamics of~$UC(i)$, which gives $\alpha_i \leq (p^* + 1)^{|UC(i)|}$.
Now~(\ref{eqn:alpha-tau-upstream})  implies that

\begin{equation}\label{alpha-upper-c<1}
\alpha = lcm\{\alpha_i : i \in [n]\}  \leq (p^* + 1)^{|CYCU|},
\end{equation}
where~$CYCU$ is as defined in Section~\ref{sec:digraphs}.  Thus,  point~(a) follows from Lemma~\ref{lem:Karp}(A3).

\medskip

For the proof of point~(b), we can use Lemma~\ref{lem:Karp}(A4).  For $d = 1$, the result follows from Proposition~\ref{prop:too-short-cycles}. For $d \geq p_*+1$, the result follows from Lemma~\ref{lem:NSC-c-fixed}(a), the second line of~\eqref{eqn:NSC-implies-alpha>1}, and~(\ref{eqn:alpha-tau-upstream}).

\medskip

For the proof of point~(c) consider $i \in [n]$. If $UC(i)$ is acyclic, then by Lemma~\ref{lem:Karp}(A1) and Proposition~\ref{prop:too-short-cycles} we may assume that $\tau_i  \leq |UC(i)| + p^* - 1 < p^* + b\ln n$ holds a.a.s.  Otherwise
$\tau_i$ is trivially bounded by the size of the state space of the internal dynamics of~$UC(i)$, which gives $\tau_i \leq (p^* + 1)^{|UC(i)|}$.
This implies that

\begin{equation}\label{tau-upper-c<1}
\tau_i \leq  \max\{(p^* + 1)^{|CYCU|}, p^* + b \ln n\} \ \mbox{a.a.s.}
\end{equation}

Now Lemma~\ref{lem:Karp}(A3) together with~(\ref{eqn:alpha-tau-upstream}) imply that $\tau = O(\ln(n))$ holds a.a.s.

On the other hand,
Lemma~\ref{lem:Karp}(A2) and Lemma~\ref{lem:path-implies-transient} imply that $\tau = \Omega(\ln(n))$ holds a.a.s. $\Box$

\subsection{Proof of Theorem~\ref{thm:c>1-upper-bd}}\label{subsec:PfThm2}

It is well known that for~$c > 1$ and any fixed~$k$, a.a.s.  property~$CYC_{[\ell]}$ will hold for some~$\ell > k$. This follows also from Lemma~\ref{lem:Karp}(B2) and the fact
that~$P(CYC_{[\ell]})$ is monotonically increasing with respect to the connection probability~$\pi(n)$.  For $p^* = th^* = 1$, the conditional probability that~$C$ is a witness of
$NSC_{[\ell]}$ given that~$C$ is a directed cycle of length~$\ell$ is~$1 - 2^{-\ell + 1}$.  In view of~\eqref{eqn:NSC-implies-alpha>1}, this implies point~(a).

\medskip
The proofs of points~(b) and~(c) are the same as for the analogous parts of Theorem~\ref{thm:c<1-upper}.
$\Box$

\subsection{Proof of Theorem~\ref{thm:c>ccrit-upper-bd}}\label{subsec:PfThm3}

By Lemma~\ref{lem:Karp}(C1) we may assume that~$D$ contains a giant component.   For each~$i$, let~$\alpha_i^-$ denote the \emph{period} of node~$i$, that is, the smallest~$T \geq 0$ such that $s_i(t) = s_i(t + T)$ for all $t \geq \tau$. Then

\begin{equation}\label{eqn:alpha-tau-DS-giant}
\begin{split}
\alpha &= lcm ( \{  \alpha_i^- : i \in DG\} \cup \{ \alpha_i :  i \notin DG \} )  \leq (p^* + 1)^{|CYCU \backslash DG|}  lcm \{\alpha_i^- : i \in DG\},\\
\tau &= \max\{\max_{i \in DG} \tau_i, \max_{i \notin DG} \tau_i \} \leq
\max \{(p^* + 1)^{|CYCU \backslash DG|}, p^* + b\ln n, \max_{i \in DG} \tau_i\} \ \mbox{a.a.s.}
\end{split}
\end{equation}

The first line of~\eqref{eqn:alpha-tau-DS-giant} follows from the same argument that we used for deriving~\eqref{alpha-upper-c<1}, and the second line is implied by ~(\ref{eqn:alpha-tau-upstream}) in analogy with~\eqref{tau-upper-c<1}.

For the proof of Theorem~\ref{thm:c>ccrit-upper-bd}, assume $p_* = p^* = p$ and~$th^* = 1$. By Lemmas~\ref{lem:ev-min-cyc}(a) and~\ref{lem:P(evmnode(g))}, a.a.s. all nodes in~$DG$ are eventually minimally cycling, which implies $\alpha_i^- = p + 1$ for every $i \in DG$.  Now point~(a) follows
from~\eqref{eqn:alpha-tau-DS-giant} and Lemma~\ref{lem:Karp}(C3).

For the proof of point~(b), assume that there exists a minimally cycling node in~$GC$. By Lemma~\ref{lem:ev-min-cyc}(b), the time~$\tau^{all}$ until all nodes in~$DG$ become minimally cycling satisfies the inequality
$\tau^{all} \leq \tau^* + L + \alpha^*(p+1)^{L+1}$, where~$\tau^* = \max_{i \notin DG} \tau_i = O(\ln(n))$ according to~\eqref{eqn:alpha-tau-DS-giant},
and~$\alpha^* \leq (p^* + 1)^{|CYCU \backslash DG|}$. Lemma~\ref{lem:Karp}(C3) shows that~$\alpha^*$ scales like a constant, and by Lemma~\ref{lem:Karp}(C5) we can assume that~$L \leq b \ln n$. Now point~(b) follows from~\eqref{eqn:alpha-tau-DS-giant} if we set~$k = b \ln (p+1)$. $\Box$

\subsection{Proof of Theorem~\ref{thm:alpha-nonpoly-lower}}\label{subsec:PfThm4}

Lemma~\ref{lem:Karp}(B1)  implies that a.a.s. the assumptions of Proposition~\ref{prop:only-1-cycle-full} are satisfied for every~$UG(i)$. Thus point~(b) follows from Proposition~\ref{prop:only-1-cycle-full}, Lemma~\ref{lem:Karp}(B3) and~\eqref{eqn:alpha-tau-upstream}.

\medskip

Similarly, Lemmas~\ref{lem:Karp}(B1),(B3) and~\ref{lem:path-implies-transient} imply point~(c).

\medskip

The proof of point~(a) requires a more substantial argument. Assume $p^* = th^* = 1$ and let $k > 0$.  We will show that

\begin{equation}\label{eqn:lcm(NSC[ell])}
\lim_{k \rightarrow \infty} \lim_{n \rightarrow \infty} P(\alpha \geq n^k) = 1.
\end{equation}

Let  $k$ be a  positive integer, let $\kappa < \beta$ be a positive constant, and let
$J = \lceil\frac{3k}{\kappa}\rceil$. Lemma~\ref{lem:Karp}(B2) assures us that for some fixed~$I$ property~$P2C(I, J, \kappa)$ will hold with probability arbitrarily close to~1 as long as~$n$ is sufficiently large.  Fix a suitable~$I$ and assume that~$P2P$ is a set that witnesses property~$P2C(I, J, \kappa)$. For simplicity of notation, assume that for each $j = q_r \ell_r \in P2P$ the node~$j$ belongs to a directed cycle of length~$j$.  By Proposition~\ref{prop:only-1-cycle-full},     if $j = q_r \ell_r \in P2P$, then the length $\alpha_{j}$ of the attractor in~$UC(j)$ is either a divisor $q_r$  or a multiple of  $\ell_{r}$. We will show that for sufficiently large $k, n$, with probability arbitrarily close to~1, we will find at least $J/3$ among the $\alpha_{j}$ that are multiples of~$\ell_{r}$. Since all $\ell_{r}$ are  pairwise  different prime numbers and satisfy $\ell_{r} > n^\kappa$, we can infer from~\eqref{eqn:alpha-tau-upstream} that

\begin{equation}\label{eqn:alpha-lb-JKL}
\alpha \geq lcm \{\alpha_{j}: \, j \in P2P\} \geq n^{\kappa J/3} \geq n^k,
\end{equation}
 and~\eqref{eqn:lcm(NSC[ell])} follows.

\medskip

For each  $j \in P2P$  consider the event~$E_j$ that   $\alpha_{j}$ is divisible by~$\ell_j$.  By Lemma~\ref{lem:Karp}(B1) we may wlog assume that~$UC(j)$ is supersimple. Thus, in particular,  no node is upstream of two distinct directed cycles, so that the sets $UC(j)$ for $j \in P2P$ are pairwise disjoint.
 This in turn implies that the events~$E_j$ are independent.   Thus the Central Limit Theorem applies, and since~$J$ can be assumed arbitrarily large, it suffices to show that each of the events~$E_j$ individually has probability larger than~$0.5$.

\begin{proposition}\label{prop:pj-divides-alphaj}
For all sufficiently large~$n$  and for each  $j = q \ell \in P2P$   the probability that $\ell$ divides $\alpha_j$ is larger than~$0.5$.
\end{proposition}

\noindent
\textbf{Proof:}  Consider a  directed cycle~$C = (i_1, i_2, \dots , i_{q\ell}, i_{q\ell+1} = i_1)$ with $q\ell$ nodes. Then $C$ can be partitioned into $\ell$ consecutive segments of length~$q$ each.  If $\ell$ does not divide the length of the attractor, then by Proposition~1 of~\cite{AJ}, at all times $t > \tau$ we must have

\begin{equation}\label{eqn:ql-attr-small}
\forall k, k' < \ell ~\, \forall m \in [q] \ s_{i_{qk + m}}(t) = s_{i_{qk' + m}}(t).
\end{equation}

Since~$q$ is odd, \eqref{eqn:ql-attr-small} would imply that
no interval $[t_0, t_0 + 2q +1)$ with~$t_0 > \tau$ can be an interval of uninterrupted firing for any node~$i \in C$. We will show that such intervals will exist with high probability, which in turn will imply that~\eqref{eqn:ql-attr-small} fails with high probability.

Note that if  $[t_0, t_0 + 2q + 1)$ is an interval of uninterrupted firings of node~$i_u$ with $u < q\ell$ and either $s_{i_u}(t_0) =  0 = s_{i_{u+1}}(t_0+1)$ or
$s_{i_u}(t_0+1) =  0 = s_{i_{u+1}}(t_0+2)$, then  $[t_0 + 1, t_0 + 2q + 2)$
will be an interval of uninterrupted firings of node~$i_{u+1}$.  Thus if~$t_0 > 0$, such an interval of uninterrupted firings will be inherited by node~$i_{u+1}$ in the next time step \emph{unless}  it gets destroyed by a firing of node $i_{u+1}$ at time~$t_0$ or~$t_0+1$.  Such a destructive firing would need to be induced by  firing input to node~$i_{u+1}$  from outside~$C$.

Let~$T$ be an arbitrary fixed positive integer.  Consider the digraph $D^-$ that is obtained from~$D$ by removing all arcs of~$C$. Since~$UC(i)$ was assumed  supersimple, for each $i \in C$ the set~$UC_{D^-}(i)$ is a tree.

For each fixed~$t>0$, let $P_i(t)$ denote the conditional probability that a given node~$i \in C$ will receive firing input at some time~$t^+ \geq t$ from outside of~$C$, given that~$C$ is a directed cycle and~$D$ is supersimple.  By symmetry, every node has the same expected properties as node~1.  Thus Lemma~\ref{lem:Karp}(B4)(B5) implies in view of   Proposition~\ref{prop:too-short-cycles} and Lemma~\ref{lem:path-implies-transient} that

\begin{equation}\label{eqn:Pt-lim-lim}
\lim_{t \rightarrow \infty} \lim_{n \rightarrow \infty} P_i(t) = 0  \qquad \mbox{and} \qquad \forall t \, \lim_{n \rightarrow \infty} P_i(t)   > 0.
\end{equation}

Thus with probability arbitrarily close to~1, there will be a time~$t^*$ such that some node~$i \in C$ receives firing input from outside~$C$ exactly $T$ more times in the interval~$[t^*, \tau]$, and no node in~$C$ receives firing input from outside~$C$  more than~$T$ times in the interval~$[t^*, \tau]$.  For different nodes~$i \in C$ the sequence of firing inputs from outside of~$C$ are independent.  Since the length of~$C$ scales like~$\Omega(n^\kappa)$ and~$T$ is fixed, for sufficiently large~$n$ it follows from~\eqref{eqn:Pt-lim-lim} that the total number of firing inputs that the nodes in~$C$ receive from outside~$C$ in the time interval~$[t^*, \infty)$ becomes a vanishingly small fraction of~$q\ell$ as $n \rightarrow \infty$.   In particular, the probability that any interval of uninterrupted firings $[t_0, t_0 + 2q + 1)$ of some node~$i$ that occurs for some~$t_0 \geq t^*$ gets subsequently destroyed becomes vanishingly small as $n \rightarrow \infty$.

Thus it suffices to show that with probability arbitrarily close to~1, for sufficiently large~$T$, there will
 be some $t_0 \geq t^*$ such that $[t_0,  t_0 + 2q + 1)$ is an interval of uninterrupted firings of some node~$i_u$.
Consider $h \notin C, i \in C$ with $<h,i> \, \in A_D$ and assume that node~$h$  fires  at time $t \geq t^*$. Then there must be a node~$v \in UC(h)$ with $s_v(0) = 0$ and a directed path $(v = h_0, h_1, \dots , h_{t} = h)$ from~$v$ to~$h$ such that
$s_{h_h}(t^-) = 0$ for all  $0 \leq t^- \leq t$. Now a straightforward inductive argument shows that if there exists also a directed path $(h_{-4q}, h_{-4q}, \dots , h_0)$
with $s_{h_\tau}(0) = 0$ iff $\tau$ is even, then we do get an interval $[t, t+ 4q + 1)$ of uninterrupted firings of node~$h$.  The latter will happen with a positive conditional probability. Since $q < I$, this conditional probability is bounded from below by a positive constant for a given choice of~$I$. Thus by
choosing~$T$ sufficiently large, we can assure that it will happen for at least one of the~$T$ firings after~$t$ with probability arbitrarily close to~1.
By Proposition~\ref{prop:ev-min-cyc2}, the interval~$[t+1, t+ 4q + 2)$ will contain an interval of uninterrupted firings of node~$i$ of length~$\geq 2q$.  The probability that this interval will subsequently be destroyed is vanishingly small, and the results follows.
$\Box$ $\Box$

\section{Open problems for future exploration}\label{sec:discussion}

Our results open a number of avenues for future research.

We investigated the average lengths of attractors~$\alpha$ and transients~$\tau$ for neuronal network models whose connectivities are Erd\H{o}s-R\'{e}nyi random digraphs with connectivity probability~$\pi(n) \approx \frac{c}{n}$. For minimal firing thresholds~$th_* =1$, our results indicate one or several phase transition for some critical values~$c_{crit} \geq 1$. Theorems~\ref{thm:c<1-upper} gives a complete characterization of the subcritical case~$c < 1$.  Theorem~\ref{thm:c>1-upper-bd}(a) indicates that at least with respect to two properties of interest, the critical value is~$c_{crit} = 1$, while Theorem~\ref{thm:c>ccrit-upper-bd} shows that with respect to two other properties the critical value satisfies the inequality   $c_{crit} \leq p^*+1$.

These results leave four open problems for the supercritical case: The first is whether the upper bound Theorem~\ref{thm:c>ccrit-upper-bd} can be improved to~$c_{crit} = 1$.  Preliminary simulation studies indicate that this is indeed the correct value, at least when $\vp = \vth = \vone$.   We include a brief summary of these in Appendix~D.
The second problem is to pinpoint the precise scaling laws for~$\tau$, as Theorem~\ref{thm:c>1-upper-bd}(c) only gives a logarithmic lower bound and Theorem~\ref{thm:c>ccrit-upper-bd}  a polynomial upper bound. Furthermore, it is not clear whether the conclusions of Theorem~\ref{thm:c>1-upper-bd}(a) and  Theorem~\ref{thm:c>ccrit-upper-bd}  continue to hold under the weaker assumption $th^* = 1$.

Theorem~\ref{thm:alpha-nonpoly-lower} covers the lower end of the ``critical window'' for $c \approx 1$.  It may be quite challenging to study the behavior in other parts of this window, such as the case~$c = 1$.  Of particular interest is the question whether for some choice of~$\pi(n)$ both the medians of~$\alpha$ and~$\tau$ can scale faster than any polynomial.

Extending the study to the case $th_* > 1$ poses another set of challenges, as it may require different techniques. The results of~\cite{JAT} imply that all conclusions of Theorems~\ref{thm:c>1-upper-bd} and Theorem~\ref{thm:c>ccrit-upper-bd} will be true for sufficiently large~$c$, but we do not know the corresponding critical values.
Theorem~\ref{thm:th*>1-trans}  suggests that a phase transition for~$\tau$ might occur at $c_{crit} = (th_*!)^{\frac{1}{th_*}}$.

Our results assume uniform distributions of initial conditions.  It may be of interest from the point of view of neuroscience applications to work out analogue results for other distributions, for example under the assumption that initially only a small fraction of nodes fire and the remaining nodes are at the end of their refractory periods. Moreover, degree distributions in actual neuronal networks may be scale-free rather than normal~\cite{Cecchi, CelegansPower} and most neuronal tissues contain neurons of several distinct types that may differ in their propensity to form connections. Thus it is of interest to study analogous problems for other types of random connectivities than Erd\H{o}s-R\'{e}nyi digraphs. Scale-free networks or inhomogeneous random digraphs as defined in~\cite{BGJ} may be particularly promising candidates for more realistic models of neuronal networks.

Finally, it is of interest to extend our study beyond lengths of attractors and transients. A particularly interesting phenomenon, know in the neuroscience context as \emph{decoherence}~\cite{ahn_terman_2010, fernandez_smith},  is the amplification of initially small Hamming distances along the transient. This is a notion of sensitive dependence on initial conditions and is recognized as a hallmark of chaotic dynamics in Boolean networks~\cite{origins}, of which our networks with $p^* = 1$ are a special case.  It would be of interest to know for which types of random networks this phenomenon is generic.

\section*{Appendix A: An example}

The following example shows that the maximum possible length of attractors and transients in our networks increases faster than~$e^{\sqrt{n}}$.   The construction modifies Examples~7 and~11 of~\cite{AJ}.

\begin{example}\label{ex:Landau}
For some function~$h(n)$ that scales like $e^{\sqrt{n\ln n } - o(1)}$ and for all~$n$ there exists a network $N_n$ on~$[n]$ with $p^* = th^* = 1$ such that $\alpha = \tau = h(n)$ for at least one initial state in~$St_{N_n}$.
\end{example}

\noindent
\textbf{Proof:}
 Consider a network on~$[n]$ with $p^* = th^* = 1$ whose connectivity~$D$  contains pairwise vertex-disjoint  directed cycle $C_\ell = (i_1^\ell, i_2^\ell, \ldots , i_{k_\ell+1}^\ell = i_1^\ell)$ for some odd $k_1, \dots , k_m > 1$ with the property that no $i_j^\ell$ is the target of an arc from outside~$C_\ell$.  Assume that in the initial state we have
$s_{i_j}^\ell = 0$ for all $\ell \in [m]$ and when $j < k_\ell$ is even, and $s_{i_j}^\ell = 1$ otherwise. Then node~$i_{k_\ell}^\ell$ will fire at all times $2t+1$ for $0 \leq t < \frac{k_\ell-1}{2}$  and will not fire at time~$ t = k_\ell$. It follows that the length~$\alpha$ of the attractor must satisfy $\alpha \geq lcm\{k_1, \dots , k_m\}$.

Now assume the network contains  another node $i^*$ such that $<i_{k_\ell}^\ell, i^*> \, \in A_D$ for all~$\ell \in [m]$. Let $t^* = lcm\{k_1, \dots , k_m\}$.  Then node~$i^*$ will receive firing input from at least one node~$i_{k_\ell}^\ell$ at all odd times $t < t^*$, but will not receive firing input at time~$t^*$ from any of the nodes $i_{k_1}^1, \ldots , i_{k_m}^m$.
Thus if $s_{i^*}(0) = 0$, then node~$i^*$ will fire at all even times $t < t^*$, and will not fire at time $t^* + 1$ unless it receives firing input from some other node of the network at time~$t^*$.

Let us complete the construction of~$D$ by adding one more directed cycle $C_0 = (i^+_0, i^+_1, i^+_0)$ of length~2, an additional arc $<i^+_0, i^*>$, and no other arcs.  In particular, all other vertices, if there are any, will be isolated.

Consider an initial state as specified above, and such that $s_{i^+_0}(0) = 0, s_{i^+_1}(0) = 1$.  Then nodes~$i^*, i^+_0$ will fire simultaneously at all even times $t < t^*$, and $s_{i^+_0}(t^*) = 1$ as $t^*$ is odd.  Thus node~$i^*$ will not receive any firing input at time~$t^*$, and it follows that $s_{i^*}(t^*+1) = 1$. Starting from time~$t^* + 2$, node~$i^*$ will fire at all odd time steps, in response to firing input from node~$i^+_0$. It follows that $\tau = t^* = lcm\{k_1, \dots , k_m\}$.

For a given~$[n]$, this construction allows us to construct networks with~$\alpha$ and~$\tau$ on the order of the maximum value of the least common multiple of odd positive integers
$k_1, \ldots , k_m$ such that $k_1 + ... + k_m \leq n -3$.  This maximum value obeys the same scaling law that was proved for Landau's function~$g(n)$ in~\cite{Landau}, and the result follows. $\Box$

\section*{Appendix B: Proofs of Lemmas~\ref{lem:Karp} and~\ref{lem:NSC-c-fixed}(b)}

Lemma~\ref{lem:Karp} mostly summarizes more or less well-known results about Erd\H{o}s-R\'{e}nyi random digraphs.  In particular, point~(C1), the second part of point~(A4), and the assertion of~(A1), (C2) that some or all $UC(i), DC(i)$ will be~$b$-small  follow directly from the results in~\cite{Karp}. Other parts of the lemma are also well-known or easily follow from well-known results, but are not explicitly stated in a single source. Point~(B2) is a custom-tool for our Theorem~\ref{thm:alpha-nonpoly-lower}.
For convenience we give here self-contained proofs of all parts of the lemma except the results of~\cite{Karp} that we mentioned earlier in this paragraph. For better flow of the exposition, our derivations are organized around four themes instead of following the order in which items appear in the text of the lemma.

\subsection*{Directed paths}

Assume  $\pi(n) \leq \frac{1 - n^{-\beta}}{n}$ for some fixed $0 < \beta <1/4$. Note that for sufficiently large~$n$ this assumption holds in part~(A) as well as in part~(B).

For each $\ell$, let $\kappa^{[\ell]}$ be the r.v. that counts the number of directed paths of length
exactly~$\ell$ in~$D$.  Then

\begin{equation}\label{eqn:E(paths-[ell])}
e^{-\ell n^{-\beta} - o(1)}n < E(\kappa^{[\ell]}) = \frac{n(n-1)\dots (n - \ell)(1-n^{-\beta})^{\ell}}{n^{\ell}} < e^{-\ell n^{-\beta}}n.
\end{equation}

 For $\ell \geq (1 + \eps)(\ln n) n^\beta$ it follows from~\eqref{eqn:E(paths-[ell])} that $P(\kappa^{[\ell]} > 0 ) \leq E(\kappa^{[\ell]}) < n^{-\eps} \rightarrow 0$. Since every directed path of length $> \ell$ contains a directed path of length~$\ell$, this implies the upper bound on~$L_{\max}$ in point~(B3).

\medskip

Consider not necessarily distinct vertices $i, j \in [n]$ and a set of arcs~$A$ in~$D$.  Let $DP(i,j, A)$ denote the event that there exists a directed path of length at least~1 from~$i$ to~$j$ that does not contain any of the arcs in~$A$, and let  $DP(i,j,A, \ell)$  denote the event that there exists such a path of length~$\ell$.

\begin{proposition}\label{prop:Pr(existspath)}
Let $0 < \beta < 1/4$, and $\pi(n) \leq \frac{1 - n^{-\beta}}{n}$.  Then for any~$A \subseteq A_D$,

\begin{equation}\label{eqn:Pr(DP(i,j))}
P(DP(i, j, A)) \leq n^{\beta - 1 + o(1)}.
\end{equation}
\end{proposition}

\noindent
\textbf{Proof:} Since the probability  $DP(i,j, A, \ell)$   of existence of a directed paths of specified length~$\ell$ from~$i$ to~$j$ is less than the expected number of such paths, we get

\begin{equation}\label{eqn:Pr(DP(i,j, ell))}
P(DP(i, j, A, \ell)) \leq \frac{(n-2)(n-3) \dots (n - \ell)(1 - n^{-\beta})^\ell }{n^\ell} <   \frac{e^{-\ell n^{-\beta}}}{n}.
\end{equation}

In particular, $P(DP(i, j, A, \ell)) < \frac{1}{n}$ for all~$\ell$. Moreover,
for  any given $\eps > 0$  and $\ell \geq n^{\beta + \eps}$ we get $P(DP(i, j, A, \ell)) < \frac{e^{-n^\eps}}{n}$. Thus from~\eqref{eqn:Pr(DP(i,j, ell))} we get

\begin{equation}\label{eqn:Pr(DP(i,j))-est}
\begin{split}
P(DP(i, j, A)) &=  \sum_{\ell = 1}^{\lfloor n^{\beta+\eps} \rfloor} P(DP(i,j,A,\ell)) +  \sum_{\ell = \lceil n^{\beta+\eps} \rceil}^{n} P(DP(i,j,A,\ell))  \\ & \leq O(e^{-n^\eps}) +  \sum_{\ell = 1}^{\lfloor n^{\beta+\eps} \rfloor} \frac{1}{n}  = O\left( \frac{n^{\beta+\eps}}{n}\right),
\end{split}
\end{equation}
and the proposition follows by considering arbitrarily small $\eps >0$. $\Box$

\bigskip

Now consider the event that an upstream or downstream component of~$D$ is supersimple.  This event will \emph{fail} to occur iff one of the following six events occurs:

\begin{itemize}
\item[(FL1)] Some node is contained in two distinct arc-disjoint directed cycles.
\item[(FL2)] Two node-disjoint directed cycles are connected by a directed path.
\item[(FL3)] There exists a directed cycle with a shortcut.  That is, there exist nodes $i \neq j$ and three pairwise arc-disjoint directed paths,  one from~$i$ to~$j$, two from~$j$ to~$i$.
\item[(FL4)] There exist distinct nodes~$i,k$ such that~$i$ belongs to a directed cycle~$C$ of~$D$ while $k \notin C$, and there are two arc-disjoint directed paths  that are both arc-disjoint from~$C$, either both from~$i$ to~$k$, or both from~$k$ to~$i$.
\item[(FL5)] Some node is downstream or upstream of two node-disjoint directed cycles.
\item[(FL6)] There exists distinct nodes~$i,j, k$ and four pairwise arc-disjoint directed paths with the following properties: A path from~$i$ to~$j$ and a path from~$j$ to~$i$ form a directed cycle~$C$. The other two paths are arc-disjoint from the directed cycle~$C$ and either go from~$i$ to~$k$ and from~$j$ to~$k$ or from~$k$ to~$i$ and from~$k$ to~$j$.
\end{itemize}

Property~(FL1) involves one node~$i$ and the existence of two arc-disjoint directed cycles that contain~$i$. If~$A$ denotes the set of arcs in the first of these direct cycles, then Proposition~\ref{prop:Pr(existspath)} implies:

\begin{equation}\label{eqn:P2adc-upper-bd}
P(FL1) \leq n P(DP(i,i, \emptyset)) P(DP(i,i,A)) \leq n \left(n^{\beta - 1 + o(1)}\right)^2 = n^{2\beta - 1 + o(1)}.
\end{equation}

\medskip

Similarly, properties (FL2)--(FL4) involve two nodes and three directed paths and properties~(FL5), (FL6) involve three nodes and four paths.  Thus the same argument that lead to~\eqref{eqn:P2adc-upper-bd} gives the inequalities.

\begin{equation}\label{eqn:P2bdc-upper-bd}
\begin{split}
P(FL2 \ \vee \ FL3 \ \vee \ FL4) &\leq n^2 \left(n^{\beta - 1 + o(1)}\right)^3 = n^{3\beta - 1 + o(1)},\\
P(FL5 \ \vee \ FL6) &\leq n^3 \left(n^{\beta - 1 + o(1)}\right)^4 = n^{4\beta - 1 + o(1)},
\end{split}
\end{equation}

Now the remaining part of point~(A1) and point~(B1) follow from our assumption on~$\beta$ and~$\pi(n)$.

\medskip

Point~(B4) and the analogue of point~(A2) for~$L_{\max}(i)$ in place of~$L_{\max}^s(i)$  are  well-known
 results about the extinction probability of the corresponding random birth process in disguise. A similar argument shows that for all fixed~$t$, with a.p.p. we have~$L_{\max}(1) \geq t$.  We can derive points~(A2) and~(B5) from these results as follows.

For a given a directed path~$Pt$ of length~$\ell$,  consider the event $Nstr(Pt)$ that~$Pt$ is not straight.  Arguing as in the proof of part~(B1) we find that

\begin{equation}\label{eqn:P(Nstr(P))}
P(Nstr(Pt)) \leq \ell^2 n^{\beta - 1 + o(1)}.
\end{equation}

Thus for $2 \leq \ell = o(n^\beta)$, the conditional probability that a sequence of nodes~$Pt$ is a straight path given that it is a path of length~$\ell$ is arbitrarily close to~$1$ for sufficiently large~$n$.  Thus points~(A2) and~(B5) follow from the analogous results about~$L_{\max}(i)$ and~$L_{\max}(1)$.

\medskip

The second part of point~(B3) could also be proved by considering extinction probabilities of birth processes and then using~\eqref{eqn:P(Nstr(P))}, but we want to give an alternative more elementary proof here.
Let $\kappa_s^{[\ell]}$ denote the number of straight paths of length~$\ell$.
 In view of~\eqref{eqn:E(paths-[ell])} this implies for any $2 \leq \ell = o(n^\beta)$ that

\begin{equation}\label{eqn:E(straight-paths-[ell])-Theta}
E(\kappa^{[\ell]}_s)  = \Theta(n).
\end{equation}

We can now use the second-moment method to derive the second part of point~(B3). We will use the following instance of Chebyshev's Inequality:

\begin{equation}\label{eqn:Chebysheff-kappa-s}
P(\kappa_s^{[\ell]}= 0) \leq P(|\kappa_s^{[\ell]} - E(\kappa_s^{[\ell]})| \geq E(\kappa_s^{[\ell]})) \leq \frac{Var(\kappa_s^{[\ell]})}{(E(\kappa_s^{[\ell]}))^2}. \end{equation}

Let us fix~$\ell \leq (1 - \eps)n^\beta$ for some $\eps > 0$.  Let $\cP$ denote the set of \emph{potential directed paths} of length~$\ell$ in~$D$, that is, the set of all sequences $Pt= (i_0, i_1, \ldots , i_\ell)$  of length~$\ell$ that consist of pairwise distinct vertices.  For each~$Pt \in \cP$,  let $\kappa_{Pt}^s$ denote the r.v. that takes the value~1 iff $Pt$ is a straight directed path in~$D$, and takes the value~$0$ otherwise. Then
$\kappa_s^{[\ell]} = \sum_{Pt \in \cP} \kappa_{Pt}^s$. Since each~$\kappa_{Pt}^s$ is an indicator variable, it follows that

\begin{equation}\label{eqn:sum-Var(kappaP)}
\sum_{Pt \in \cP} Var(\kappa_{Pt}^s) \leq E(\kappa_s^{[\ell]}).
\end{equation}

In order to be able to use~\eqref{eqn:Chebysheff-kappa-s}, we need to estimate the covariances. Let $\cP^2_s$ denote the set of all pairs $(Pt_1, Pt_2) \in \cP^2$ such that $Pt_1 \neq Pt_2$ and $Pt_1, Pt_2$ contain exactly one nonempty maximal stretch of consecutive common potential arcs. For the purpose of this definition we will consider a common vertex as a ``common stretch'' of length~0. Notice that if  $Pt_1, Pt_2$ are two different elements of~$\cP$ such that
$(Pt_1, Pt_2) \notin    \cP^2_s$, then there are two possibilities: Either~$Pt_1$ and~$Pt_2$ are arc-disjoint, which implies that~$\kappa_{Pt_1}^s$ and $\kappa_{Pt_2}^s$ are independent and have covariance~0. Or $Pt_1, Pt_2$ contain more than one such stretch, which implies that~$\kappa_{Pt_1}^s$ and~$\kappa_{Pt_1}^s$ cannot simultaneously take the value~1 and thus have negative  covariance.  Since the variables $\kappa^s_{Pt}$ take nonnegative values, together with~\eqref{eqn:sum-Var(kappaP)},  this observation implies the following estimate.

\begin{equation}\label{eqn:Var(kappa-s-ell)-est2}
\begin{split}
Var(\kappa_s^{[\ell]})\  &\leq \  E(\kappa_s^{[\ell]}) \ + \ \sum_{(Pt_1, Pt_2) \in \cP^2_s} Cov(\kappa_{Pt_1}^s, \kappa_{Pt_2}^s)\\
&\leq \ E(\kappa_s^{[\ell]})\  + \ \sum_{(Pt_1, Pt_2) \in \cP^2_s} E(\kappa_{Pt_1}^s \kappa_{Pt_2}^s).
\end{split}
\end{equation}

Let us estimate the right-hand side of~\eqref{eqn:Var(kappa-s-ell)-est2}.

Fix~$Pt_1 \in \cP$ and consider the set $\{Pt_2 \in \cP:\ (Pt_1, Pt_2) \in \cP^2_s\}$. For each $Pt_2$ in this set, there are $\frac{\ell(\ell +1)}{2} - 1$ possible locations of the endpoints of the common stretch in~$Pt_1$. These locations in~$Pt_1$ determine the length of the common stretch.  Thus given these endpoints, there are at most~$\ell$ possibilities of choosing the corresponding locations in~$Pt_2$, which gives a rough estimate of fewer than $\ell^3$ possible locations of the common stretch in the pair~$(Pt_1, Pt_2)$. For a fixed choice of these locations and a common stretch with~$k < \ell + 1$ vertices (and hence $k-1$ arcs), we have

\begin{equation}\label{eqn:number-of-Pt2s}
(n-\ell-1)\dots (n - 2\ell + k - 1) < n^{\ell - k + 1}
\end{equation}
choices for the corresponding potential paths $Pt_2$ (as they are determined by choosing the remaining~$\ell - k + 1$ vertices in~$Pt_2$), while
for each such choice the conditional probability $P(\kappa_{Pt_2}^s = 1| \kappa_{Pt_1}^s = 1)$ satisfies

\begin{equation}\label{eqn:P(kappaP2|kappaP1)}
P(\kappa_{Pt_2}^s = 1| \kappa_{Pt_1}^s = 1) \leq \frac{(1 - n^{-\beta})}{n^{\ell - k + 1}} < \frac{1}{n^{\ell - k + 1}}.
\end{equation}

It follows from~\eqref{eqn:number-of-Pt2s} and~\eqref{eqn:P(kappaP2|kappaP1)} that for any given $Pt_1 \in \cP$ we have

\begin{equation}\label{eqn:Var(kappa-s-ell)-est4}
\sum_{\{Pt_2 \in \cP:\ (Pt_1, Pt_2) \in \cP^2_s\}} E(\kappa_{Pt_1}^s \kappa_{Pt_2}^s) \leq P(\kappa_{Pt_1}^s = 1) \ell^3.
\end{equation}

Note that  $E(\kappa^{[\ell]}_s) = \sum_{P \in \cP} P(\kappa_P^s = 1)$. Thus~\eqref{eqn:Var(kappa-s-ell)-est4} in turn implies, in view of~\eqref{eqn:E(paths-[ell])} and our choices of $\beta$ and~$\ell$, that

\begin{equation*}
\sum_{Pt_1 \in \cP}\quad  \sum_{\{Pt_2 \in \cP:\ (Pt_1, Pt_2) \in \cP^2_s\}} E(\kappa_{Pt_1}^s \kappa_{Pt_2}^s) \leq   \ell^3 E(\kappa^{[\ell]}_s)  =
 O(n^{1 + 3\beta}).
\end{equation*}

From~\eqref{eqn:E(straight-paths-[ell])-Theta}  and~\eqref{eqn:Var(kappa-s-ell)-est2}      we can infer that

\begin{equation*}
Var(\kappa_s^{[\ell]}) = O(n) + O(n^{1+ 3\beta}) = O(n^{1+ 3\beta}) = o(E(\kappa_s^{[\ell]})^2),
\end{equation*}
and~\eqref{eqn:Chebysheff-kappa-s} implies  for
  $\ell \leq (1-\eps) n^{\beta}$ that

\begin{equation*}
\lim_{n \rightarrow \infty} P(\kappa_s^{[\ell]} > 0) = 1,
\end{equation*}
which proves the second part of point~(B3).

\medskip

\subsection*{Directed cycles}

We call a sequence $C = (i_1, i_2, \ldots ,  i_k, i_1)$ such that $i_1, i_2, \ldots , i_k$ are pairwise distinct a \emph{potential directed cycle.}
Let~$\cC$ denote the set of all  potential directed cycles.
For~$C \in \cC$, define a r.v. $\xi_C$ by letting $\xi_C = 1$ if all potential arcs in~$C$ are in~$A_D$,  and $\xi_C = 0$
otherwise.

If~$\pi(n) = \frac{c}{n}$,  the expected value of the r.v. $\xi_C$ is
$E(\xi_C)= \frac{c^k}{n^k}$.

Since each directed cycle in~$D$ of length~$\ell$ is represented by~$\ell$ potential directed cycles that differ by cyclic shifts, the number $\xi$ of directed cycles  in~$D$ is given by~$\xi = \sum_{C \in \cC} \frac{\xi_C}{|C|}$.

Let $\eta$ be the number of nodes that belong to directed cycles, and let~$\eta^{[\ell]}$ denote the number of nodes that belong to directed cycles of length~$\ell$.
Assume $c < 1$.
Then

\begin{equation}\label{eqn:E(nodes-in-cycles)-k}
\begin{split}
E(\eta^{[\ell]}) &\leq  \ell\frac{n(n-1) \cdots (n-\ell+1)}{n^\ell} \frac{c^\ell}{\ell} \leq {c^\ell},\\
E(\eta) & \leq \sum_{\ell = 2}^\infty E(\eta^{[\ell]}) \leq \sum_{\ell = 2}^\infty c^\ell < \infty.
\end{split}
\end{equation}

This implies point~(A3).

\medskip

For fixed $\ell, \ell' \geq 2$ let $\cC^{[\ell]}$ denote the set of all potential directed cycles of length~$\ell$, and let $\cC^{\ell, \ell'}$ denote the set of all ordered pairs  $(C, C')$ such that $C \in \cC^{[\ell]}, C' \in \cC^{[\ell']}, C \neq C'$, and
$C, C'$ have at least one arc  in common.   Similarly, let $\cC^{\ell, \ell'}_+$ denote the set of all ordered pairs $(C, C')$ such that $C \in \cC^{[\ell]}, C' \in \cC^{[\ell']}, C \neq C'$, and
$C, C'$ have at least one vertex in common. Consider  additional r.v.s

\begin{equation}\label{eqn:define-eta2-ell}
\xi^{[\ell, \ell']} = \sum_{  (C, C') \in \cC^{\ell, \ell'}} \xi_C\xi_{C'}, \qquad
\xi^{[\ell, \ell']}_+ = \sum_{(C, C') \in \cC_+^{\ell, \ell'}} \xi_C\xi_{C'},
\qquad
\xi_2 = \sum_{\ell = 2}^n \sum_{\ell' = 2}^n  \xi^{[\ell, \ell']}_+.
\end{equation}

\begin{proposition}\label{prop:E(xiellell')}
(a) Let $\ell, \ell' \geq 2$ be fixed and assume $\pi(n) \leq \frac{c}{n}$ for some fixed~$c$.  Then

\begin{equation}\label{eqn:E(xiellell')}
E\left(\xi^{[\ell, \ell']}\right) = O\left(n^{-1}\right).
\end{equation}

\smallskip

\noindent
(b) If $\pi(n) \leq \frac{1 - n^{-\beta}}{n}$ for some fixed~$\beta$ with $0 < \beta < 1/4$, then

\begin{equation}\label{eqn:E(xi[2])}
E\left(\xi_2\right) = n^{2\beta - 1 + o(1)}.
\end{equation}
\end{proposition}

\noindent
\textbf{Proof of part~(a):}  Part~(b) can be derived from Proposition~\ref{prop:Pr(existspath)} in exactly the same way as~\eqref{eqn:P2adc-upper-bd} above. In view of part~(b),  part~(a) is of interest for our purposes only for $c \geq 1$, but we prove that it holds in general.

\smallskip

We need the following version of Exercise~6.44 in~\cite{bookchapter}.

\begin{proposition}\label{prop:common-arc-vertices}
Let $C, C' \in \cC$ with $|C| = \ell, |C'| = \ell'$, and assume that $C, C'$ have $M$ vertices and $m$ potential arcs in common.
If $M > 0$ then either $M > m$ or $C = C'$.
\end{proposition}

\noindent
\textbf{Proof:}  Let $S$ be the set of all sources~$i$ for which there exists a common potential arc $<i, j>$ in both $C$ and $C'$, and let~$T$  be the set of
targets~$j$ of all common potential arcs $<i, j>$ of $C$ and $C'$. Then $m = |S| =|T|$. If $S = T$, then we must have $C = C'$ or $m = 0$.  If not, then there exists $j \in T \backslash S$, and the result follows from the observation that $M \geq |S \cup T|$. $\Box$

\bigskip

For each $C = (i_1, i_2, \dots , i_\ell, i_1) \in \cC^{[\ell]}$,  and each subset $I \subseteq \{i_1, \ldots , i_\ell\}$, let
$\cC^{\ell, \ell'}_I(C)$ denote the set of all $C' \in \cC^{[\ell']} \backslash \{C\}$ such that $I$ is the intersection of the sets of vertices of $C$ and~$C'$.  Furthermore, let $\cC^{\ell, \ell'}(C) = \bigcup \{\cC^{\ell, \ell'}_I(C): \ I \subseteq \{i_1, \ldots, i_\ell\} \ \& \ |I|  \geq 2 \}$, and for $2 \leq M \leq \ell$ let
$\cC^{\ell, \ell'}(C, M) = \bigcup \{\cC^{\ell, \ell'}_I(C): \ I \subseteq \{i_1, \ldots, i_\ell\} \ \& \ |I| = M\}$.  Then

\begin{equation}\label{eqn:E(xi2-ell-ell')-est1}
E(\xi^{[\ell, \ell']}) =  \sum_{C \in \cC^{[\ell]}} \ \sum_{C' \in \cC^{\ell, \ell'}(C)} E(\xi_C\xi_{C'}) =
 \sum_{   C \in \cC^{[\ell]}   } \ \sum_{M=2}^\ell \ \sum_{C' \in \cC^{\ell, \ell'}(C, M)} E(\xi_C\xi_{C'}).
\end{equation}

Now consider $C$ and $C' \in \cC^{\ell, \ell'}(C, M)$.   The r.v. $\xi_{C} \xi_{C'}$ takes the value~1 if, and only if, all potential arcs of both~$C$ and~$C'$ are in~$A_D$.  Thus if $m$ is the number of common arcs of~$C$ and~$C'$, then

\begin{equation}\label{eqn:E(xiC-timesetaC'-upper-m}
E(\xi_{C}\xi_{C'}) = P(\xi_C\xi_{C'} = 1) < \frac{c^{\ell + \ell' - m}}{n^{\ell + \ell' - m}}.
\end{equation}

By Proposition~\ref{prop:common-arc-vertices}, if $M > 0$ and $C \neq C'$, it follows from~\eqref{eqn:E(xiC-timesetaC'-upper-m} that we have

\begin{equation}\label{eqn:E(xiC-timesetaC'-upper-M}
E(\xi_{C}\xi_{C'}) = P(\xi_{C}\xi_{C'} = 1) < \frac{c^{\ell + \ell' - M + 1}}{n^{\ell + \ell' - M + 1}}.
\end{equation}

Moreover, note that there are  $\binom{\ell}{M} < \ell!$ possible locations for the set of common vertices in~$C$ and
$\binom{\ell'}{M} < (\ell')!$ possible locations for the set of common vertices in each of~$C'$.
Thus  equations~(\ref{eqn:E(xi2-ell-ell')-est1}) and~(\ref{eqn:E(xiC-timesetaC'-upper-M}) imply:

\begin{equation}\label{eqn:E(xi2-ell-ell')-est2}
\begin{split}
E(\xi^{[\ell, \ell']}) &<
 \sum_{C \in \cC^{[\ell]}} \sum_{M=2}^\ell \ell!(\ell')! \frac{c^{\ell+ \ell' -M + 1}}{n^{\ell + \ell' -M + 1}}  \\
& <  \ell!\ell'!\sum_{M=2}^\ell \frac{c^{\ell + \ell' -  M + 1}n^{\ell + \ell' - M}}{n^{\ell + \ell'- M + 1}}
=  \ell!\ell'!\sum_{M=2}^\ell
\frac{c^{\ell + \ell' -  M + 1}}{n},
\end{split}
\end{equation}
and point~(a) follows.
$\Box$

\bigskip

Now we can derive the remaining part of point~(A4) of the lemma and show that~$CYC_{[\ell]}$ has a.p.p. First note that we need to prove the result only for the case $\pi(n) = \frac{c}{n}$ with $1 \geq c > 0$ as  $P(CYC_{[\ell]})$  increases monotonically with respect to~$c$.

\medskip

For each~$\ell$ choose a subset $\cU^{[\ell]} \subset \cC^{[\ell]}$ of unique representations of potential directed cycles that contains for each $C \in \cC^{[\ell]}$ exactly one cyclic permutation of~$C$. Now fix $\ell \geq 2$ and  consider a r.v.

\begin{equation}\label{eqn:define-eta-ell-first}
\xi^{[\ell]} = \sum_{C \in \cU^{[\ell]}} \xi_C.
\end{equation}

Then $CYC_{[\ell]}$ holds iff $\xi^{[\ell]} > 0$. We have

\begin{equation}\label{eqn:E(eta-ell)}
E(\xi^{[\ell]}) =     \frac{n(n-1) \cdots (n-\ell+1)c^\ell}{ \ell  n^\ell}  \approx
 \frac{c^\ell}{\ell },
\end{equation}
where the approximation becomes arbitrarily good for fixed~$\ell$ as $n \rightarrow \infty$.

Now consider the  r.v.

\begin{equation}
\xi^{[\ell]}_2 = \sum_{C, C' \in  \cU^{[\ell]}, C \neq C'} \xi_C\xi_{C'}
\end{equation}

Note that if $\xi^{[\ell]} = k > 1$, then $\xi^{[\ell]}_2 = k(k-1) \geq 2$, and it follows that

\begin{equation}\label{eqn:E(eta-ell)-E(eta2-ell)}
E(\xi^{[\ell]}) \leq P(\xi^{[\ell]} > 0) + E(\xi^{[\ell]}_2).
\end{equation}

From Proposition~\ref{prop:E(xiellell')}(a) we obtain

\begin{equation}\label{eqn:E(eta2-ell)-est2}
\begin{split}
E(\xi^{[\ell]}_2) &=  \sum_{C \in  \cU^{[\ell]}} \quad \sum_{C' \in  \cU^{[\ell]}, (C,C') \notin C^{\ell, \ell}, C' \neq C }   E(\xi_C \xi_{C'})
+ \sum_{C \in \cU^{[\ell]}} \quad \sum_{ C'\in \cU^{[\ell]},  (C,C') \in C^{\ell, \ell}}   E(\xi_C \xi_{C'})   \\
&<  \frac{n(n-1) \dots (n-2\ell+1)c^{2\ell}}{\ell^2  n^{2\ell}}
+ \sum_{C \in \cU^{[\ell]}} \quad \sum_{ C'\in \cU^{[\ell]},  (C,C') \in C^{\ell, \ell}} E(\xi_C \xi_{C'}) \  \\
&< \frac{c^{2\ell}}{\ell^2 } + \sum_{C \in \cC^{[\ell]}} \quad \sum_{  C'\in \cU^{[\ell]},  (C,C') \in C^{\ell, \ell} }  E(\xi_C \xi_{C'}) = \frac{c^{2\ell}}{\ell^2 } +
O\left(n^{-1}\right).
\end{split}
\end{equation}

Now it follows from~\eqref{eqn:E(eta-ell)}
and~\eqref{eqn:E(eta-ell)-E(eta2-ell)} that

\begin{equation}\label{eqn:P(eta[ell])-estimate}
\lim_{n \rightarrow \infty} P(\xi^{[\ell]} > 0) \geq  \frac{c^{\ell}}{ \ell } -  \frac{c^{2\ell}}{  \ell^2 }  > 0,
\end{equation}
where the last inequality follows from our working assumption that $c \leq 1$.  This implies that~$CYC_{[\ell]}$ has a.p.p. $\Box$

\subsection*{The proof of point~(B2)}

Assume~$\pi(n) = \frac{1 - n^\beta}{n}$ for some fixed~$0 < \beta < \frac{1}{4}$.
We  will need the following sharper version of~\eqref{eqn:P(eta[ell])-estimate}.

\begin{proposition}\label{prop:lambda-exists}
 If  $0 < \lambda < \beta$, then  for all  $2 \leq \ell \leq  n^\lambda$:

\begin{equation}\label{eqn:P(eta[ell])-estimate-L}
 P(\xi^{[\ell]} > 0) = \frac{1- o(1)}{\ell} .
\end{equation}
\end{proposition}

\noindent
\textbf{Proof:} Similarly to~\eqref{eqn:E(paths-[ell])} we get  the estimate

\begin{equation}\label{eqn:E(eta-ell)-beta1}
\frac{e^{-\ell n^{-\beta} + o(1)}}{\ell} < E(\xi^{[\ell]}) =    \frac{n(n-1) \cdots (n-\ell+1)(1 - n^{-\beta})^\ell}{\ell n^\ell}  < \frac{e^{-\ell n^{-\beta}}}{\ell}.
\end{equation}

It follows that every fixed $w < 1$, all $\ell$ as in the assumption, and sufficiently large~$n$

\begin{equation}\label{eqn:E(eta-ell)-beta-final}
\frac{w}{\ell} < E(\xi^{[\ell]}) < \frac{1}{\ell}.
\end{equation}

Now we can proceed as in the proof at the end of the previous subsection.
Instead of part~(a) we can now use part~(b) of Proposition~\ref{prop:E(xiellell')}. This gives us the following modification of~\eqref{eqn:E(eta2-ell)-est2}

\begin{equation*}\label{eqn:E(eta2-ell)-est-beta}
\begin{split}
E(\xi^{[\ell]}_2) &<  \sum_{C \in \cC^{[\ell]}} \quad \sum_{C' \in  \cU^{[\ell]}, (C,C') \notin C^{\ell, \ell}_+, C' \neq C }    E(\xi_C \xi_{C'})\ + \sum_{C \in \cC^{[\ell]}} \quad \sum_{  C'\in \cU^{[\ell]},  (C,C') \in C^{\ell, \ell}_+}  E(\xi_C \xi_{C'})  \ \\
& < \frac{n(n-1) \cdots (n-2\ell+1)c^{2\ell}}{ \ell^2  n^{2\ell}}  + \sum_{C \in \cC^{[\ell]}} \quad \sum_{  C'\in \cU^{[\ell]},  (C,C') \in C^{\ell, \ell}_+}  E(\xi_C \xi_{C'})       \\
&< \frac{1-o(1)}{\ell^2} + E(\xi_2) =  \frac{1-o(1)}{\ell^2}  + n^{2\beta - 1 + o(1)},
\end{split}
\end{equation*}
and the result follows from~\eqref{eqn:E(eta-ell)-E(eta2-ell)}  and~\eqref{eqn:E(eta-ell)-beta-final}. $\Box$

\bigskip

Now  consider, for $\ell$ as above, r.v.'s $\zeta^{[\ell]} = \min \{1,  \xi^{[\ell]} \}$ that take the value~1 if $CYC_{[\ell]}$ holds and the value~0 otherwise.
Let us define  $K = \lfloor n^\kappa \rfloor-1, L = \lceil n^\lambda\rceil-1$ for some fixed  $0 < \kappa < \lambda < \beta$.

For integers $q,  K, L$ with $1 \leq q  < K  < L \leq \frac{n}{q}$, let $Primes(K, L)$ denote the set of prime numbers~$\ell$ with $K < \ell < L$.
Consider the r.v.s

\begin{equation}\label{eqn:define-eta-prime}
\zeta_{q, K, L} = \sum_{\ell \in Primes(K, L)} \zeta^{[q\ell]}.
\end{equation}

 By Proposition~\ref{prop:lambda-exists}, for each fixed~$q$ as above we have

\begin{equation}\label{eqn:E(zeta)-q-K-L-lim}
\lim_{n \rightarrow \infty} E(\zeta_{q, K, L}) =  \lim_{n \rightarrow \infty} \frac{1 - o(1)}{q} \sum_{\ell \in Primes(K, L)}
 \frac{1}{\ell}.
\end{equation}

According to  refinements of Euler's famous result as derived in~\cite{meissel, mertens},  we have

\begin{equation}\label{eqn:Euler}
\lim_{L \rightarrow \infty} \left(  \left(\sum_{\ell \in Primes(2, L)}  \frac{1}{\ell} \right)  -   \ln(\ln(L))   \right) = M \approx 0.2615,
\end{equation}
where~$M$ is known as the \emph{Meissel-Mertens constant.}

In view of our choice of $K , L$, we conclude that for every fixed $q$ we will have

\begin{equation}\label{eqn:E(zeta)-q-limn}
\lim_{n \rightarrow \infty}  E(\zeta_{q, K,L}) = \lim_{n \rightarrow \infty} \frac{\ln(\ln(L)) - \ln(\ln(K ))}{q}
=  \frac{\ln  \lambda - \ln  \kappa}{q}.
\end{equation}

Note that since  $0 < \kappa < \lambda < 1$, the number $\ln \lambda  - \ln  \kappa$ is positive.

Now let us fix $\kappa, \lambda, K, L$ as above and define, for each integer  $I > 1$,
 a r.v. $\zeta_I$ as the sum of the r.v.s $\zeta_{q, K, L}$ over
all prime numbers~$2 < q < I$.  Then, again by (\ref{eqn:Euler}) and~\eqref{eqn:E(zeta)-q-limn},

\begin{equation}\label{eqn:E(zeta)-limn-inf}
\lim_{I \rightarrow \infty} \lim_{n \rightarrow \infty} E(\zeta_I) =
\lim_{I \rightarrow \infty} \lim_{n \rightarrow \infty} \sum_{q \in Primes(2,I)} E(\zeta_{q, K, L}) =\infty.
\end{equation}

The r.v. $\zeta_I$ is a sum of indicator variables, and it follows that

\begin{equation}\label{eqn:Var(zeta-I())}
Var(\zeta_I) \leq  E(\zeta_I) \ + \ 2\sum_{q, q' \in Primes(2, I)} \quad \sum_{\ell, \ell' \in Primes(K, L)} Cov(\zeta^{[q\ell]},  \zeta^{[q'\ell']}),
\end{equation}
where the sum of covariances is taken over all pairs $(q, \ell) \neq (q', \ell')$.

Since the r.v.s~$\zeta^{[q\ell]}$ take values in~$\{0, 1\}$ we have

\begin{equation}\label{eqn:Cov-sum-zeta}
\begin{split}
 Cov(\zeta^{[q\ell]},  \zeta^{[q'\ell']}) &< P(\zeta^{[q\ell]}\zeta^{[q'\ell']} > 0) = P(\xi^{[q\ell]}\xi^{[q'\ell']} > 0)\\
 &< \sum_{C \in \cC^{[q\ell]}} \sum_{C' \in \cC^{[q'\ell']}}  E(\xi_C \xi_{C'}) < E(\xi_2).
\end{split}
\end{equation}

It follows  that the sum of covariances in~\eqref{eqn:Var(zeta-I())} is bounded from above by~$E(\xi_2)$, which goes to zero as $n \rightarrow \infty$ in view of
Proposition~\ref{prop:E(xiellell')}(b).

Thus   $\lim_{n \rightarrow \infty}  (Var(\zeta_I) -  E(\zeta_I)) = 0$   and it follows from Chebyshev's Inequality and~\eqref{eqn:E(zeta)-limn-inf} that for every fixed positive integer~$J$ and $\eps > 0$ we can choose a positive integer~$I$ such that

\begin{equation}\label{eqn:J-large}
\lim_{n \rightarrow \infty} P(\zeta_I > J) > 1 - \eps.
\end{equation}

This completes the proof of part~(B2).

\subsection*{The supercritical case}

Assume $c > 1$ and consider a randomly chosen~$D$. By~\cite{Karp}, we may assume  that~$D$ contains a giant strongly connected component~$GC$  such that $DG$ has size of size approximately $\varrho(c) n$. By symmetry, we may assume without loss of generality that $[n] \backslash  DG = [m]$.  Then there will be no arcs from~$DG$ to $[m]$ in~$D$, and the arcs from~$[m]$ into~$DG$ are irrelevant for properties~(C2)--(C4). Thus we are considering properties of a random digraph on~$m$, with $\pi(m) = \frac{c}{n} \approx \frac{(1- \varrho(c))c}{m}$, where
$m \approx (1 - \varrho(c))n$ approaches infinity as~$n$ does.  As long as  $(1 - \varrho(c)) c < 1$, this allows us to deduce property~(C2) from property~(A1), property~(C3) from property~(A3), and property~(C4) from property~ (A2).

Thus it  suffices to show that $(1 - \varrho) c < 1$, where $\varrho$ is the unique root of the equation

\begin{equation}\label{eqn:varrho}
e^{-c\varrho} = 1 - \varrho
\end{equation}
in the interval~$(0,1)$.

For $c > 1$ the function $f(x) = 1 - x - e^{-cx}$ takes positive values on the interval $[0, \varrho)$ and negative values on the interval~$(\varrho, 1]$. Moreover,
 $f(1 - \frac{1}{c}) = \frac{1}{c} - e^{-c + 1} > 0$ iff $g(c) = ce^{-c + 1} < 1$. Since $\frac{dg}{dc} = (1 - c)e^{-c + 1} < 0$ and   $g(1) = 1$, we conclude that
 $1 - \frac{1}{c} < \varrho$, which is equivalent to $(1 - \varrho) c < 1$.

\medskip

Finally, let us prove~(C5). It is well known that the analogous result holds for undirected random graphs. That is, if $c > 1$ and edges are randomly and independently drawn with probability~$\frac{c}{n}$, then for some constant $a = a(c) > 0$, with probability approaching~1 as $n \rightarrow \infty$,  the giant, and in fact all components will have diameter $\leq a\ln n$ (see, for example, \cite{Bollobas}, \cite{Durrett}, or~\cite{Luczak}).  The standard proof of this result can be adapted to the directed case with minor modifications, but we will give an alternative argument here that reduces~(C5) to the result for undirected graphs.

Consider the following procedure for producing a random graph~$G$ on~$[n]$: First draw a random digraph~$D$ with~$\pi(n) = \frac{c}{n}$, where $c > 1$. Now use a breadth-first search to produce a random graph by  recursively adding vertices to a list~$V = \bigcup_{\ell, m} V_\ell^m $ and edges to the set~$E$. Each of the sets $V_0^m$ will be of the form $V_0^m = \{i^m\}$ and $V_\ell^m$ will consist of all vertices that can be reached in~$D$ from $i^m$
by a directed path of length~$\ell$, but not by a shorter directed path.   Initially, $V := V_1^0 := \{1\}$.  Given a current list $V$, consider the last subset $V_\ell^m$ that was added to~$V$.  Let $V^m_{\ell + 1}$ be the set of all
$j \in [n] \backslash V$ such that there exists $i \in V^m_\ell$ with $<i, j> \, \in A_D$. If $V^n_\ell \neq \emptyset$, add all the aforementioned edges $\{i, j\}$ to $E$, as well as edges $\{j_1, j_2\} \subset V^m_{\ell+1}$ for which $j_1 < j_2$ and $<j_1, j_2> \, \in A_D$, and add $V^m_{\ell+1}$ to~$V$.  If
 $V^n_\ell = \emptyset$ and $[n] \backslash V \neq \emptyset$,  let $V_0^{m+1}  = \{i_{m+1}\},$   where~$i_{m+1}$ is the lowest-numbered vertex outside of~$V$.  Repeat until $V = [n]$.

Since it doesn't matter for the probability of obtaining a specific~$G$ by this procedure whether~$D$ is chosen at the outset or membership of $<i, j>, <j,i>$ in~$A_D$ is determined  at the time of  deciding whether $\{i, j\} \in E$,  this procedure produces a random Erd\H{o}s-R\'{e}nyi graph~$G$ with edge probability $\frac{c}{n}$. The same construction can be repeated multiple times for a randomly chosen permutation of the nodes, and since node~$1$ has positive probability of ending up in the giant strongly connected component of~$D$, we may in the remainder of this argument assume that node~$1$ is in~$GC_D$.

The set $V^1 = \bigcup_{\ell} V^1_\ell$ is the connected component of~$1$ in~$G$ and is equal to~$DC(1)$ in~$D$. Thus with probability that can be chosen arbitrarily close to~$1$, by the result for undirected graphs,  every node~$j \in DC(1)$ can be reached in~$D$ from~$1$ by a directed path of length~$\leq a \ln n$.  The dual construction shows that,  also with probability that can be chosen arbitrarily close to~$1$, node~$1$ can be reached in~$D$  from every every node~$i \in UC(1)$  by a directed path of length~$\leq a \ln n$. Thus property~(C5) follows for  any $b \geq 2a$.  $\Box$

\subsection*{The proof of Lemma~\ref{lem:NSC-c-fixed}(b)}

Consider~$\ell \geq k$ such that~$\ell$  is a multiple of~$p+1$. Define random variables $\zeta_C$ for $C  = (i_0, i_1, \dots , i_\ell = i_0) \in \cC$  that take the value~1 iff all potential arcs in~$C$ are in~$A_D$, the potential directed cycle $C$ is composed exclusively of nodes~$j$ with $p_j = p$ and $th_j = 1$, and $s_{i_j}(0) = j + r \ mod \ (p+1)$ for some fixed $r \in \{0, \ldots, p\}$ and  all $i_j \in C$, and takes the value~0 otherwise.
Then

\begin{equation}\label{eqn:E(zetaC-p)}
E(\zeta_C) =   \frac{c^\ell}{n^\ell  (p+1)^{\ell-1} (p^*-p_* + 1)^\ell (th^*)^\ell} = \frac{ (p+1) c^\ell}{m^\ell n^\ell},
\end{equation}
where $m =  (p+1)(p^*-p_* + 1) th^* \geq 2$.

Consider the r.v. $\zeta^{[\ell]} = \sum_{C \in  \cU^{[\ell]}} \zeta_C$.
Note that if~$\zeta_C = 1$, then some cyclic shift of~$C$ witnesses~$NSC_{[\ell]}(p)$.  Thus the inequality $\zeta^{[\ell]} > 0$ implies $NSC_{[\ell]}(p)$.  We have

\begin{equation}\label{eqn:E(zetap-ell)}
E(\zeta^{[\ell]}) =    \frac{n(n-1) \cdots (n-\ell+1)  (p+1) c^\ell}{ \ell  m^\ell n^\ell}  \approx \frac{ (p+1)  c^\ell}{ \ell  m^\ell},
\end{equation}
where the approximation becomes arbitrarily good for fixed~$\ell$ as $n \rightarrow \infty$.

If $C, C'$ are  arc-disjoint, then $\zeta_C$ and $\zeta_{C'}$ are independent  and hence $Cov(\zeta_C, \zeta_{C'}) = 0$.  If
$C, C'$ have a common arc, that is, if $(C_1, C_2) \in C^{\ell, \ell}$, then  $Cov(\zeta_C, \zeta_{C'}) \leq P(\zeta_C = 1)P(\zeta_{C'} = 1) \leq E(\xi_C\xi_{C'})$.
Since each variable~$\zeta_C$ is an indicator variable, it follows from~\eqref{eqn:E(zetap-ell)} and Proposition~\ref{prop:E(xiellell')}(a) that

\begin{equation}\label{eqn:Var-zetaC-estim}
Var(\zeta^{[\ell]}) = \sum_{C \in \cC^{[\ell]}} Var(\zeta_C) +  O\left(n^{-1}\right) \leq \sum_{C \in \cC^{[\ell]}} E(\zeta_C) +  O\left(n^{-1}\right)
\leq \frac{ (p+1)  c^\ell}{m^\ell} +  O\left(n^{-1}\right).
\end{equation}

Thus

\begin{equation}\label{eqn:var(zeta)/E(zeta)-squared}
\frac{Var(\zeta^{[\ell]})}{E(\zeta^{[\ell]})^2} \leq \frac{\ell^2 m^\ell}{ (p+1)  c^\ell} + o(1).
\end{equation}

If $c > m$ and~$\ell$ sufficiently large, by substituting~\eqref{eqn:var(zeta)/E(zeta)-squared} in Chebyshev's Inequality we can deduce that
$P(NSC_{[\ell]}(p))$ is arbitrarily close to~1.   $\Box$

\section*{Appendix C: Proof of Theorem~\ref{thm:th*>1-trans}}\label{subsec:PfThm5}

\subsection*{Necessary and sufficient conditions for long transients if $th_* > 1$}\label{subsec:transients-for-thm4}

In Subsection~\ref{subsec:transients} we investigated the existence of transients for the case~$th_* = 1$; here we develop counterparts of these results for the case $th_* > 1$  that are needed for the proof of Theorem~\ref{thm:th*>1-trans}.

By Proposition~\ref{prop:too-short-cycles}, if~$D$ is acyclic, the existence of long transients in~$N$ requires the existence of sufficiently long directed paths in~$D$. If $th_* > 1$, then certain types of trees are required.  To see this, let us consider again a node $i$ and assume $s_i(t^*) = 0$ for some~$t^* > 0$. A \emph{branched history of~$s_i(t^*) = 0$} is a function
$F: TR \rightarrow [n]$, where $TR$ is a set of finite sequences of positive integers of lengths $\leq t^*$ that is closed under subsequences such that $F(\emptyset) = i$ and

\begin{itemize}
\item[(TR0)] $s_{F(\sigma)}(t^* - |\sigma|) = 0$ for all $0 \leq |\sigma| \leq t^*$. In particular, $s_{F(\emptyset)} (t^*) = s_i(t^*) = 0$.
\item[(TR1)] If $k \in [th_{F(\sigma)}]$ for some $\sigma$ with~$|\sigma| < t^*$, then $\sigma^\frown k \in  TR$ and $<F(\sigma^\frown k), F(\sigma)> \, \in A_D$.
\end{itemize}

The definition of the dynamics of~$N$ requires that every firing has at least one branched history.  The structure $(TR, \subseteq)$ is a tree with root $\emptyset$ whose leaves are the sequences of length~$t^*$.  Let $Tr(t^*)$ be the tree of all such sequences of length $\leq t^*$ that take values in the set~$[th_*]$.   Note that $Tr(t^*)$ must be a subtree of the domain $TR$ of any branched history of any firing at time~$t^*$.  If~$D$ is acyclic, then~(TR1) implies that any branched firing history~$F: TR \rightarrow [n]$ is an injection.  This proves the following:

\begin{proposition}\label{prop:branched-hist-acyc}
Assume $D$ is acyclic and let~$d$ be the largest  integer such that there is an embedding of $Tr(d)$ into $D$. Then
$s_i(d+1) \geq 1$ for all~$i \in [n]$.  In particular, $\tau \leq d + p^*$.
\end{proposition}

Proposition~\ref{prop:branched-hist-acyc} will allow us to derive upper bounds for ~$\tau$ when $th_* > 1$. The following construction provides a tool for deriving  lower bounds.

\begin{definition}\label{def:witness}
 Let $d$ be a positive integer and consider the set~$Tr = Tr(d)$ of all finite sequences~$\sigma = (\sigma_1, \sigma_2, \ldots , \sigma_{|\sigma|})$ of length $0 \leq |\sigma| \leq d$ such that
$\sigma_i \in [th_*]$ for all~$i$.

For fixed $D$ and $\vs = \vs(0)$ call $W \subseteq [n]$ a \emph{witness of depth~$d$} if there exists a bijection $F: Tr \rightarrow W$ such that

\begin{itemize}
\item[(i)] $p_i = p_*$ and $th_i = th_*$ for all $i \in W$.
\item[(ii)] $s_{F(\sigma)} = s_{F(\sigma^\frown k)} - 1 \mod (p_*+1)$ for all $k \in [th_*]$  and $\sigma \in Tr$ with $|\sigma| < d$. Moreover, $s_{F(\sigma)} = 0$ when $|\sigma| = d$.
\item[(iii)] $\{j: \ <j, F(\sigma)> \, \in A_D\} = \{F(\sigma^\frown k): \ k \in [th_*]\}$  for all $\sigma \in Tr$ with $|\sigma| < d$;
\item[(iv)] $\{j: \ <j, F(\sigma)> \, \in A_D\} = \emptyset$  for all $\sigma \in Tr$ with $|\sigma| = d$;
\item[(v)] $F(\sigma^\frown k) < F(\sigma^\frown (k+1))$  for all $k \in [th_*]$ and $\sigma \in Tr$ with $|\sigma| < d$.
\item[(vi)] $A_D$ has no arcs with endpoints in the range of~$F$ other than the ones required by~(iii).
\end{itemize}

\end{definition}

The above definition of a witness is stronger and more technical than is strictly needed for the proof of the next proposition, but the additional conditions will make it easier to estimate the probability of existence of witnesses.

\begin{proposition}\label{prop:witness-implies-transient}
Let $p_*, p^*, th_*, th^*$ be arbitrary. Suppose there exists a witness of depth~$d$. Then $\tau \geq d +1$.
\end{proposition}

\noindent
\textbf{Proof:} Condition~(v) ensures  that if $W$ is a witness, then there exists a unique bijection that satisfies conditions (iii)--(vi).  In particular, $F(\emptyset)$ is uniquely determined by~$W$ itself, and conditions~(iii), (iv) imply that $UC(F(\emptyset)) = W$ and the subdigraph of $D$ that is induced by $W$ is acyclic. Thus by
Proposition~\ref{prop:too-short-cycles}, all nodes in $UC(F(\emptyset)) = W$ will eventually stop firing in the trajectory of~$\vs(0)$. Thus the length of the trajectories will be $\geq t^*+1$, where~$t^*$ is the last time at which any node in~$W$ fires.   On the other hand, a simple inductive argument shows that conditions~(i)--(iii) imply that
$s_{F(\sigma)}(t) = 0$ when $|\sigma| = d - t$. In particular, $s_{F(\emptyset)}(d) = 0$. This shows that  $\tau \geq d + 1$. $\Box$

\subsection*{Proof of Theorem~\ref{thm:th*>1-trans}}

Let $d$ be a positive integer, and let $Tr = Tr(d)$ be as in Definition~\ref{def:witness}. Then $(Tr, \subseteq)$ forms a tree with root $\emptyset$ such that

\begin{equation}\label{eqn:TR-size}
|Tr| = \frac{th_*^{d+1} - 1}{th_* - 1}.
\end{equation}

\medskip

Now let $F: Tr \rightarrow [n]$ be a randomly chosen injection, and let $W$ denote its range. Let us estimate the probability that $F$ satisfies conditions~(i)--(vi) of Definition~\ref{def:witness}, which imply that~$W$ is a witness.

Since the vectors $\vp$ and $\vth$ are drawn from a uniform distribution subject to constraints on minimal and maximal values, the event that $F$ satisfies~(i) has probability

\begin{equation}\label{eqn:P(i)-min-trans}
P((i)) = \left(\frac{1}{(p^*-p_*+ 1)(th^*-th_*+ 1)}\right)^{|Tr|} \leq 1.
\end{equation}

Similarly, the conditional probability that $F$ satisfies~(ii) given that it satisfies~(i) is

\begin{equation}\label{eqn:P(ii|i)-min-trans}
P((ii)|(i)) = \left(\frac{1}{p_*+ 1}\right)^{|Tr|}.
\end{equation}

The probability that~(iv) holds for a given $\sigma$ with $|\sigma| = d$ is approximately $e^{-c}$; more generally, since we need to exclude both incoming and outgoing arcs with endpoints in~$W$, the probability that~(vi) holds as well
is approximately

\begin{equation}\label{eqn:P(vi)-min-trans}
P((iv) \ \& \ (vi)) \approx e^{-2c|Tr|}.
\end{equation}

Since~(iii) requires one outgoing arc for every~$F(\sigma)$ except $\sigma = \emptyset$, the probability that all arcs required by~(iii) are in~$A_D$ is equal to

\begin{equation}\label{eqn:P(iii)-min-trans}
P((iii)) =  \left(\frac{c}{n}\right)^{|Tr| - 1}.
\end{equation}

Since the relevant events except~(i) and~(ii) are independent, the conditional property that $F$ satisfies (i)--(vi) if it satisfies~(v) can be calculated as

\begin{equation}\label{eqn:P(i-vi)-min-trans}
P((i)-(vi)|(v))  = \frac{K^{|Tr|}n^{1 - |Tr|}}{c},
\end{equation}
where

\begin{equation}\label{eqn:def-K}
K =  \frac{ce^{-2c}}{(p^*-p_*+ 1)(th^*-th_*+ 1)(p_*+ 1)} < 1.
\end{equation}

In the remainder of this argument it will be convenient to sometimes write $m$ for  $th_*$ in order to reduce clutter.

Now let $\xi$ be the random variable that counts  the number of witnesses.   Since each witness uniquely determines  its embedding~$F$, the value of~$\xi$ is equal to the number of bijections~$F$ for which~(i)-(vi) hold.  There are $n(n-1)\dots (n - |Tr| + 1)$ injections of $Tr$ into~$[n]$. The fraction of them that satisfy condition~(v) is $(m!)^{m^d-|Tr|}$. Thus the expected value of~$\xi$ satisfies

\begin{equation}\label{eqn:E(Tr-witnesses)}
E(\xi) = \frac{n(n-1)\dots (n - |Tr| + 1)(K)^{|Tr|}n}{n^{|Tr|}(m!)^{|Tr|-m^d}c} = \frac{(K + o(1))^{|Tr|}n}{(m!)^{|Tr|-m^d}c},
\end{equation}
as long as $|Tr| = o(n)$.

Define

\begin{equation}\label{eqn:def-K2}
K_1 =  \frac{(K)^{\frac{m}{m-1}}}{(m!)^{\frac{1}{m-1}}} < 1.
\end{equation}

By~(\ref{eqn:TR-size}),  Equation~(\ref{eqn:E(Tr-witnesses)}) translates into

\begin{equation}\label{eqn:E(Tr-witnesses)-m>1}
E(\xi) = \frac{1}{c} \frac{(K + o(1))^{\frac{m^{d+1} - 1}{m - 1}}n}{(m!)^{\frac{m^d-1}{m-1}}} = \frac{n}{c} \left(\frac{m!}{K+o(1)}\right)^{\frac{1}{m-1}}(K_1 + o(1))^{m^d}.
\end{equation}

If  $d = \log_m (b \ln  n)$, where $b$ is any positive constant, then~(\ref{eqn:E(Tr-witnesses)-m>1}) simplifies to

\begin{equation}\label{eqn:E(Tr-witnesses)-m}
E(\xi) =  \frac{1}{c} n^{1 + b\ln (K_1 + o(1))},
\end{equation}
and it follows that

\begin{equation}\label{eqn:E(Tr-witnesses)-m-infty}
th_* = m > 1 \ \& \  b < \frac{1}{-\ln K_1} \ \& \ d  = \log_m (b \ln n) \Rightarrow \lim_{n \rightarrow \infty} E(\xi) = \infty.
\end{equation}

By~(\ref{eqn:def-K}) and~(\ref{eqn:def-K2}), the bounds on~$b$ for which~(\ref{eqn:E(Tr-witnesses)-m-infty}) works  depend only on $c, th_*, th^*, p_*, p^*$.   Thus point~(a)  of the theorem will follow from Proposition~\ref{prop:witness-implies-transient} if we can show that

\begin{equation}\label{eqn:E-inf-2nd}
\lim_{n \rightarrow \infty} E(\xi) = \infty \Rightarrow \lim_{n \rightarrow \infty} P(\xi = 0) = 0.
\end{equation}

This follows from fairly general properties of the r.v.~$\xi$.  Notice that $\xi = \sum_F \xi_F$, where $F$ ranges over all injections of~$Tr$ into $[n]$ that satisfy~(v), and $\xi_F$ takes the value~$1$ if $F$ satisfies~(i)--(vi) and takes the value~$0$ otherwise.  Thus we have~$Var(\xi_F) \leq E(\xi_F)$ for all~$F$.
Moreover, condition~(vi) ensures that if $F \neq G$, then either their ranges are disjoint and the r.v.s  $\xi_F, \xi_G$ are independent, or their ranges overlap and the
r.v.s  $\xi_F, \xi_G$ cannot simultaneously take the value~$1$.  Thus we have $Cov(F, G) \leq 0$ whenever $F \neq G$, and it follows that $Var(\xi) \leq E(\xi)$.  By Chebyshev's Inequality,

\begin{equation}\label{eqn:Chebysheff}
P(\xi = 0) \leq P(|\xi - E(\xi)| \geq E(\xi)) \leq \frac{Var(\xi)}{(E(\xi))^2} \leq \frac{1}{E(\xi)},
\end{equation}
and~(\ref{eqn:E-inf-2nd}) follows.

\medskip

For the proof of~(b), note that $m = th_* > 1$  implies that every firing must have a (not necessarily injective) branched history and at least   one node  must fire at some time $t^*$ with $t^* + p_* > \tau$. Thus it suffices to estimate an upper bound for the probability
that there is an function~$F: Tr \rightarrow [n]$ that satisfies
\begin{itemize}
\item[(iii-)] $\{j: \ <j, F(\sigma)> \, \in A_D\} \subseteq \{F(\sigma^\frown k): \ k \in [th_*]\}$  for all $\sigma \in Tr$ with $|\sigma| < d$;
\item[(v)] $F(\sigma^\frown k) < F(\sigma^\frown (k+1))$  for all $\sigma \in Tr$ with $|\sigma| < d$ and $k < th_*$.
\end{itemize}

Now consider the r.v. $\eta$  that represents the number of such functions.  We have $n$ choices for $F(\emptyset)$. The function~$F$ does not need to be injective, but for each of the $|Tr| - m^d$ nodes~$\sigma \in Tr$ that are not leaves we need to pick~$m$ pairwise distinct successors that satisfy~(v). In analogy  with~\eqref{eqn:E(Tr-witnesses)} we get the estimate

\begin{equation}\label{eqn:E(Tr-weak-witnesses)}
E(\eta) = \frac{n(n(n-1)\dots (n - m + 1))^{|Tr|-m^d}c^{|Tr|}}{n^{|Tr|}(m!)^{|Tr|-m^d}c} = \frac{n(c + o(1))^{|Tr|}}{(m!)^{|Tr|-m^d}c},
\end{equation}
as long as $|Tr| = o(n)$. As before, this translates into

\begin{equation}\label{eqn:E(Tr-weak-witnesses)-m>1}
E(\xi) = \frac{(c + o(1))^{\frac{m^{d+1} - 1}{m - 1}}n}{(m!)^{\frac{m^d-1}{m-1}} c} = \left(\frac{m!}{c}\right)^{\frac{1}{m-1}} (K_2 + o(1))^{m^d} \frac{n}{c},
\end{equation}
where

\begin{equation}\label{eqn:def-K3}
K_2 =  \frac{c^{\frac{m}{m-1}}}{(m!)^{\frac{1}{m-1}}}.
\end{equation}

If  $d = \log_m (b \ln  n)$, where $b$ is any positive constant,  $c < (m!)^{\frac{1}{m}}$, and $K_2 \neq 1$, then~(\ref{eqn:E(Tr-weak-witnesses)-m>1}) simplifies to

\begin{equation}
E(\xi) =  \frac{1}{c} n^{1 + b\ln (K_2 + o(1))}.
\end{equation}

If $c < (m!)^{\frac{1}{m}}$, then  $\ln K_2 < 0$,  and we conclude

\begin{equation}\label{eqn:E(Tr-weak-witnesses)-m-infty}
th_* = m > 1 \ \& \ b > \frac{1}{-\ln K_2} \ \& \ d  = \log_m (b \ln n) \Rightarrow \lim_{n \rightarrow \infty} E(\xi) = 0.
\end{equation}

Point~(b) now follows from the inequality $P(\xi > 0) \leq E(\xi)$. $\Box$

\section*{Appendix D: Preliminary simulation results}

We numerically explored the lengths~$\alpha$ of the attractors and~$\tau$ of the transients for the case $\vec{p} = \vec{th} =\vec{1}$ near the critical value $\pi(n) = 1/n$. For this purpose we wrote a \textsc{MatLab} program that works as follows:

\begin{itemize}
\item \textbf{Input:}
\begin{itemize}
\item The number $n$ of nodes in the network.
\item The vector of connection probabilities $\pi(n) = \frac{c}{n}$ to be explored.
\item The number $N$ of repetitions per parameter setting.
\end{itemize}
\item\textbf{Output:}
\begin{itemize}
\item The observed values of~$\alpha$ and~$\tau$.
\end{itemize}
\end{itemize}

We ran this program for $n \in \{100, 200, 400, 800, 1600, 3200\}$ and $\pi(n) = \frac{c}{n}$ with $c \in \{0.8, 0.82, \ldots, 1.48, 1.5\}$. For each parameter setting we did $2,000$ repetitions, where each repetition used a new randomly drawn digraph $D$ and a new initial condition $\vec{s}(0)$ that was randomly drawn from the uniform distribution.

Figure~1 summarizes the findings that are most relevant to this paper.  It shows the median, maximum, and $99.9$th percentiles  of the observed values of~$\alpha$ and~$\tau$ for each parameter setting that we explored. The results roughly confirm the theoretical predictions of our theorems.

The peaks of the medians, shown in Panels~(A) and~(D) respectively, appear to occur for~$c > 1$ and move closer to~1 as the number of nodes~$n$ increases.  This adds some empirical evidence to our conjecture that at least for  $p_* = p^* = 1$  the result of Theorem~\ref{thm:c>ccrit-upper-bd} can be improved by weakening the assumption~$c > p+1$ to~$c > 1$. This also may indicate that the most complex dynamics occurs near the upper end of the critical window for emergence of the giant connected component.

However, for no parameter setting did we observe median attractor lengths that exceed the median lengths of transients, as would be predicted by Theorem~\ref{thm:alpha-nonpoly-lower}.  The most likely explanation appears to be that while the median of~$\alpha$ is predicted by the theorem to grow superpolynomially, the effect may be observable only for much larger~$n$ than we were able to explore numerically.   Panels~(B-C) and~(E-F)  give some supporting evidence for this explanation. For larger $n$ the maximum observed values of~$\alpha$ are of a higher order of magnitude than the maximum observed values for~$\tau$.  Maximum values are sensitive to outliers and rare events.  It may well be the case that these ``rare'' events become common for much larger~$n$.  To explore this possibility, we computed  $99.9$th percentile  (interpreted as the mean of the second and third largest) values of $\alpha$ and $\tau$. Panels (E-F) show that  the values for~$\alpha$ are significantly larger than  the ones of $\tau$,  but the effect shows up only for larger $n$.

\begin{figure}
 \centering
 	 \includegraphics[width=6in,height=5in]{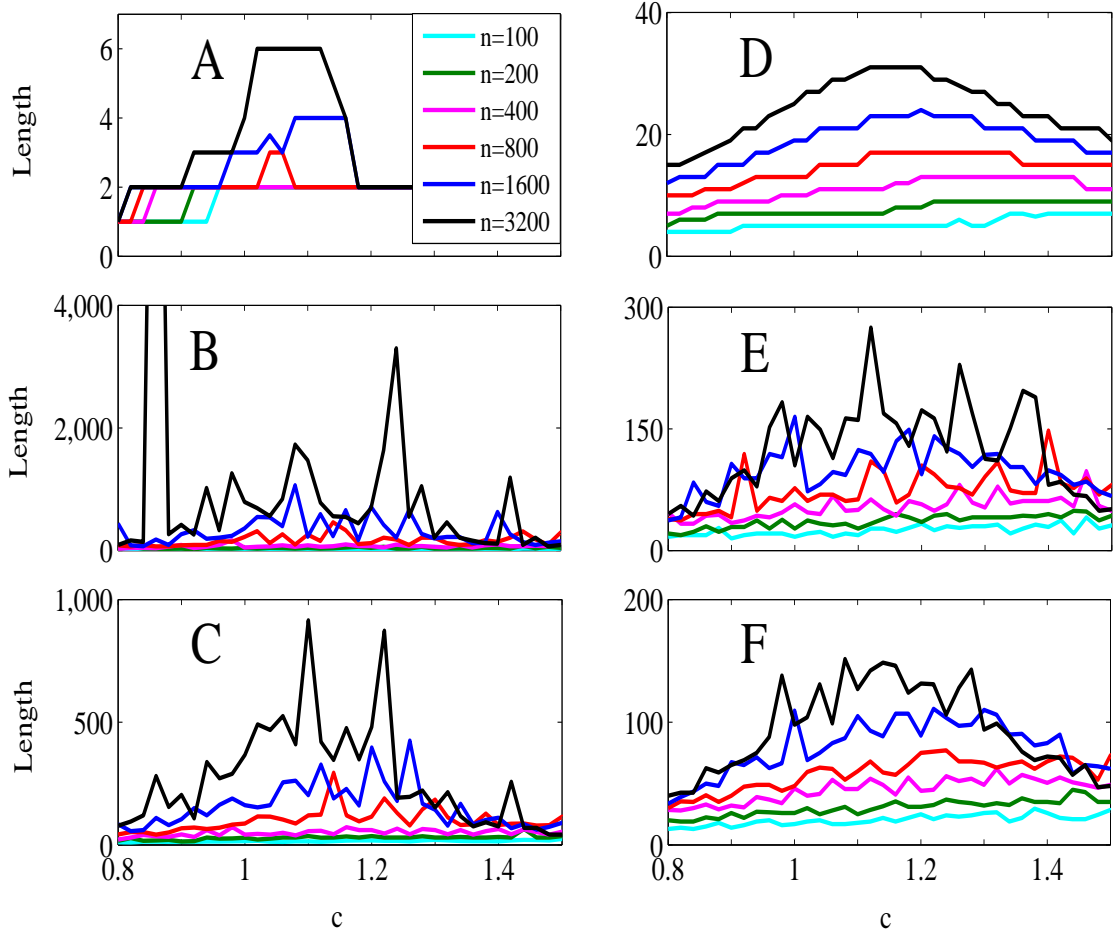}
 \caption{\label{fig_att_trans} Lengths of the attractors and the transients for $\vec{p} = \vec{th} =\vec{1}$.  (A-C) Median, maximum, and $99.9$th percentile of $\alpha$.  (D-F) Median, maximum, and  $99.9$th percentile of~$\tau.$ }
\end{figure}

\end{document}